\documentclass{amsart} 
\usepackage{graphicx}
\usepackage{varioref}

\newtheorem{theorem}{Theorem}[section] 
\newtheorem{lemma}[theorem]{Lemma}

\newtheorem{corollary}[theorem]{Corollary}

\theoremstyle{definition}
\newtheorem{definition}[theorem]{Definition}
\newtheorem{example}[theorem]{Example} 
\newtheorem{problem}[theorem]{Problem} 

\theoremstyle{remark}
\newtheorem{remark}[theorem]{Remark}
\newtheorem{ack}{Acknowledgements}

\renewcommand{\bar}{\overline}
\newcommand{\A}{\ensuremath{\mathcal A}}

\newcommand{\D}{\ensuremath{\mathcal D}}

\newcommand{\X}{\ensuremath{\mathcal X}}
\newcommand{\C}{\ensuremath{\mathbb C}}
\renewcommand{\L}{\ensuremath{\mathsf{L}}}

\newcommand{\G}{\ensuremath{\mathcal G}}
\newcommand{\sG}{\ensuremath{\mathfrak{G}}}
\renewcommand{\NG}{\ensuremath{\mathcal {NG}}}
\newcommand{\NP}{\ensuremath{\mathcal {NP}}}
\newcommand{\pV}{\ensuremath{\overline{V}}}
\newcommand{\R}{\ensuremath{\mathcal R}}
\newcommand{\V}{\ensuremath{\Sigma}}
\renewcommand{\Re}{\ensuremath{\mathbb R}}

\newcommand{\Cir}{\ensuremath{\mathcal C}}
\newcommand{\Z}{\ensuremath{\mathbb Z}}
\renewcommand{\P}{\ensuremath{\mathbb P}}
\newcommand{\bl}{\ensuremath{{\boldsymbol{\lambda}}}}
\renewcommand{\l}{\ensuremath{\lambda}} 

\newcommand{\we}{\ensuremath{\wedge}}
\newcommand{\be}{\ensuremath{{\boldsymbol{\eta}}}}
\newcommand{\e}{\ensuremath{\eta}}

\newcommand{\bx}{\ensuremath{{\boldsymbol{\xi}}}}
\newcommand{\x}{\ensuremath{\xi}}
\newcommand{\bm}{\ensuremath{{\boldsymbol{\mu}}}}
\newcommand{\m}{\ensuremath{\mu}}
\newcommand{\bn}{\ensuremath{{\boldsymbol{\nu}}}}
\newcommand{\n}{\ensuremath{\nu}}
\newcommand{\codim}{\operatorname{codim}}
\newcommand{\charac}{\operatorname{char}}
\newcommand{\supp}{\operatorname{supp}}

\newcommand{\Ann}{\operatorname{Ann}}
\newcommand{\depth}{\operatorname{depth}}
\newcommand{\rank}{\operatorname{rank}}
\newcommand{\Hom}{\operatorname{Hom}}

\begin{document}
\title{The line geometry of resonance varieties} 
\author{Michael Falk}

\address{Department of
Mathematics and Statistics\\ Northern Arizona University \\ Flagstaff, AZ 86011-5717}
\email{michael.falk@@nau.edu} 

\keywords{resonance variety, characteristic variety, line complex, arrangement, matroid, Orlik-Solomon algebra, local system cohomology}
\subjclass[2000]{14N20, 16S99, 14J26, 05B35}

\begin{abstract}
Let $\R^1(\A,R)$ be the degree-one resonance variety over a field $R$ of a hyperplane arrangement \A.
We give a geometric description of $\R^1(\A,R)$ in terms of projective line complexes. 
The projective image of $\R^1(\A,R)$ is a union of ruled varieties, parametrized by neighborly partitions of subarrangements of $\A$. The underlying line complexes are intersections of special Schubert varieties, easily described in terms of the corresponding partition. We generalize the definition and decomposition of $\R^1(\A,R)$ to arbitrary commutative rings, and point out the anomalies that arise. In general the decomposition is parametrized by neighborly graphs, which need not induce neighborly partitions of subarrangements of \A.

We use this approach to show that the resonance variety of the Hessian arrangement over a field of characteristic three has a nonlinear component, a cubic threefold with interesting line structure. This answers a question of A.~Suciu. We show that Suciu's deleted $B_3$ arrangement has resonance components over $\Z_2$ that intersect nontrivially. We also exhibit resonant weights over $\Z_4$ supported on the deleted $B_3,$  which has no neighborly partitions. 
The modular resonant weights on the deleted $B_3$ exponentiate to points on the complex torus which lie on, and determine, the translated 1-torus in the first characteristic variety. 
\end{abstract}
\maketitle

\begin{section}{Resonance and characteristic varieties}

Arising out of the study of local system cohomology and fundamental groups, characteristic and resonance varieties of complex hyperplane arrangements have become the object of much of the current research in the field. The study of resonance varieties in particular has led to surprising connections with other areas of mathematics: Kac-Moody algebras, Latin squares and loops, nets, special pencils of plane curves, homological algebra, ${\mathfrak sl}_2$ representations, critical points,  and Fuchsian differential equations \cite{LY,Y8, F11,EPY,SchSuc,SchSuc2, CV}. This paper adds to the list, relating resonance varieties to the theory of projective line complexes.

Rank-one local systems $\mathcal{L}={\mathcal L}_{\bf t}$ on a space $X,$ with coefficients in a field $R,$ are parametrized by $\Hom(H_1(X,\Z),R^*)\cong (R^*)^n,$ where $n$ in the first betti number of $X$. The $R$-characteristic varieties of $X$ are defined by $$\V^d_k(X,R)=\{{\bf t}\in (R^*)^n \ | \ \dim_R H^d(X,{\mathcal L}_{\bf t})\geq k\}.$$ We will restrict our attention to the case where $X$ is the complement in $\C^\ell$ of the union of a finite set  \A\ of linear hyperplanes, though some of the results cited below hold for arbitrary quasi-projective varieties.
Characteristic varieties originated, in case $d=1$  and $R=\C$, in work of Libgober on Alexander invariants \cite{L4} and were studied in that setting by several authors \cite{L5,Hironaka,CS3}. The stratification $\V^1_k(X,\C), \, k\geq 0$ of $(\C^*)^n$ determines the first betti numbers of finite abelian covers of $X$, among which is  the Milnor fiber of the non-isolated singularity of $\bigcup\A$ at the origin.
Characteristic varieties over finite fields were considered in \cite{Matei,MatSuc1,MatSuc2}. These determine numerical invariants of $\pi_1(M)$, including $p$-torsion ($p\not = \charac(R)$) in the first homology of finite abelian covers of $X$. For any  quasiprojective variety $X$,  the components of $\V^d_1(X,\C)$ are cosets of subtori of $(\C^*)^n$ by elements of $(S^1)^n$; for positive-dimensional components of  $\V^1_1(X,\C)$ the translating elements have finite order \cite{Arap}.

Rank-one complex local systems also play a role in the theory of generalized hypergeometric functions \cite{OT2}. Here one is interested in the top cohomology $H^\ell(X,{\mathcal L}_{\bf t})$.  For ${\bf t}$ satisfying some genericity conditions, $H^*(X,{\mathcal L}_{\bf t})$ can be computed from $H^*(X,\C)$ \cite{ESV,STV,L5}. The latter is treated as a cochain complex with differential $d_\bl$ given by left multiplication by $\bl\in H^1(X,\C)\cong \C^n,$ with $\exp(2\pi i \bl)={\bf t}$. This motivated the definition of resonance varieties in \cite{F7}, as the support loci in $H^1(X,\C)$ for the cohomology of $H^*(X,\C)$  relative to $d_\bl.$ 
Specifically, $$\R^d(X,R)=\{\bl \in H^1(X,R) \ | \ H^d(H^*(X,R),d_\bl)\not =0\}.$$ 
Again we have a related stratification $\{\R^d_k(X,R), \ | \ k\geq 0\}$ of $\C^n$ for each $d\geq 0$, given by $R^d_k(X,R)=\{\bl \in H^1(X,R) \ | \ \dim H^d(H^*(X,R),d_\bl)\geq k\}$. 
The case $R=\Z_N$ yields information about complex local systems: if $\bf t$ is a rational point on $\V^1_1(X,\C)$, then $t=\exp(2 \pi i \bl/N)$ where $\bl$ is a $\Z_N$-resonant weight \cite{CO}. (There are issues with the interpretation of the latter statement when $N$ is not prime - see Remark \ref{norank}.)
 The resonance variety $\R^d_k(X,\C)$ coincides with the tangent cone at the identity to $\V^d_k(X,\C)$ \cite{CS3,LY,CO,L6}. Thus the \C-resonance varieties are unions of linear subspaces, by the result of \cite{Arap} cited above, and determine the components of the corresponding characteristic varieties passing through 1. The tangent cone result fails in positive characteristic  - see \cite[Example 10.7]{Suc2}.

Since they are defined in terms of the cohomology ring, the resonance varieties of the complement of a complex hyperplane arrangement depend only on the underlying matroid of the arrangement \cite{OS1}. It is not known whether the characteristic varieties are so determined - this is a major open question. The problem is to identify by some combinatorial means components that do not contain the identity. Examples of such translated components are somewhat rare. The first to be found is an isolated point of $\V^1_2(X,\C),$ for $X$ the complement of the non-Fano arrangement \cite{MatSuc1}. The first positive-dimensional example is a translated 1-torus in $\V^1_1(X,\C),$ for $X$ the complement of the deleted $B_3$ arrangement \cite{Suc}. Several other positive-dimensional examples have been found \cite{Suc,Suc2}, including an infinite family \cite{Coh2}. We will see that the non-Fano and deleted $B_3$ examples arise from, or at least reflect, the incidence combinatorics of the underlying matroids; we suspect the same is true for all the other known examples as well. The same incidence structure gives rise to resonant weights in positive characteristics - see Section \ref{sect-examples}.

Little is understood about $\R^d(X,\C)$ for $d>1.$ It is shown in \cite[Theorem 4.1(b)]{EPY} that resonance ``propogates," that is, $\R^d(X,R)\subseteq \R^{d+1}(X,R)$ . There are resonant local systems which are not resonant in degree one \cite{CV}.

Degree-one resonance varieties over a field $R$ of characteristic zero
can be calculated directly \cite{F7}, or can be understood in terms of the Vinberg-Kac classification of Cartan matrices \cite{LY}. The latter approach gives an alternate proof that the components of $\R^1(X,R)$ are linear, and in addition shows that they intersect trivially. For arrangements of projective lines in $\P^2(\C),$ $\C$-resonant weights give rise, via the theory of ruled surfaces, to pencils of  curves among whose singular elements are unions of the lines of $\A$ \cite{LY,F11}. These determine  partitions of the arrangement with very special properties \cite{LY}. Special types of degree-one resonance varieties are related to nets, loops,  and the group law on the nonsingular cubic \cite{Y8}, and to $K(\pi,1)$ arrangements \cite{F11}. 
There are descriptions of $\R^1(X,\C)$ in terms of the linearized Alexander matrix \cite{CS3}, or via a chain complex derived from stratified Morse Theory \cite{CO}. Among all of the various approaches to $\C$-resonance varieties extant, only the direct method of \cite{F7} can be extended to fields of positive characteristic.

\medskip
It is in this rich context that we study the geometry of degree-one resonance varieties of complements of arrangements of hyperplanes. Our main interest is in resonance over fields of positive characteristic. We give a decomposition of the resonance variety into combinatorial pieces and show that, projectively, each of these pieces is the ruled variety corresponding to an intersection of (special) Schubert varieties in special position in the Grassmannian of lines in projective space. The definition of resonance varieties is extended to arbitrary commutative rings and the combinatorial decomposition is shown to hold in this generality. We show that the resonance variety of the Hessian arrangement, over an algebraically-closed field of characteristic three, has nonlinear components, irreducible cubic threefolds with interesting geometry. We show that the resonance variety of the deleted $B_3$ arrangement, over a field of characteristic two, has (linear) components which have nontrivial intersection. As noted above, neither of these phenomena can occur over fields of characteristic zero. The deleted $B_3$ arrangement also has essential resonant weights over $\Z_4,$ with full support, which do not yield the neighborly partitions that characterize resonance over fields. We relate the latter phenomena to the positive-dimensional translated component in the first characteristic variety of the deleted $B_3.$

Here is a more detailed outline of the main results of the paper. In Section \ref{resonance} we describe our main objects of study and extend the definition of degree-one resonant weights and resonance varieties to arbitrary commutative rings. Let $\A=\{H_1,\ldots, H_n\}$ be an arrangement of hyperplanes, with underlying matroid ${\sG},$ and let $A=A({\sG},R)$ denote the Orlik-Solomon algebra of $\sG$ over a commutative ring $R$. Then $A$ is a graded-commutative $R$-algebra (with no torsion) generated by degree-one elements $a_1,\ldots, a_n.$
An element $a=\sum_{i=1}^n \l_i a_i$ is called {\em resonant} if there is an element $a' =\sum_{i=1}^n \e_i a_i\in A^1$ such that $a\we a'=0$ in $A^2$, with the determinantal rank of $\begin{bmatrix}\bl | \be\end{bmatrix}$ equal to two. The coefficient vector $(\l_1,\ldots,\l_n)$ is called a {\em resonant weight.} The collection of all resonant weights forms the degree-one  {\em resonance variety} $\R^1(A,R)$. This agrees with the definition above when $R$ is a field.
If $R$ is an integral domain, we define $\R_k^1(A)=\R_k^1(A\otimes F)\cap R^n,$ where $F$ is the field of quotients of $R$. Then $\R_1^1(A,R)=\R^1(A,R).$ If $R$ is not an integral domain, there seems to be no definition of $\R^1_k(A)$ for $k>1$ such that the resulting set consists of resonant weights in the sense we have adopted.

In Section \ref{decomp} we derive a combinatorial description of resonant pairs of weights in degree one. This results in a decomposition of $\R^1(A)$ into subvarieties $V^1(\Gamma,R)$ determined by graphs $\Gamma$ with vertex set $[n]:=\{1,\ldots,n\}$. The graphs  are necessarily {\em \sG-neighborly} (Def. \ref{neighborlydef}). If $R$ is an integral domain, the decomposition of $\R^1(A)$ extends to $\R_k^1(A),$ for all $k\geq 1.$ In this case neighborly graphs induce partitions of submatroids of $\sG$ which are neighborly in the sense of \cite{F7}.

In Section \ref{sect-lines} we describe the geometry of $V^1(\Gamma,R)$  in terms of projective line geometry. Here, for simplicity, we work over an algebraically closed field.
Let  $\pV(\Gamma,R)$ be the projective image of $V^1(\Gamma,R)$. Let $K=K(\Gamma,R)$ be the kernel in $R^n$ of the incidence matrix with columns indexed by points $i\in[n]$ and rows by nontrivial rank-two flats $X$ of $\sG$ which are not contained in cliques of $\Gamma$.
For each maximal clique $S$ in $\Gamma$, let $D_S=\{\bx\in K \ | \x_i=0 \ \text{for all} \ i\in S\}$. Let $\D=\D(\Gamma)$ be the resulting subspace arrangement, and $\overline{\D}$ the corresponding projective subspace arrangement. For $D\in \D$, let $\L_D$ denote the set of lines in $\P^{n-1}$ which meet the projective image of $D$; then $\L_D$ is a subvariety of the Grassmannian $\G(2,n)$. Let $\L(\D)=\bigcap_{D\in \D} \L_D.$ Then $\pV(\Gamma,R)$ is precisely the union of the lines in $\L(\D)$. We call the subspaces in $\overline{\D}$ the {\em directrices} of  $\pV(\Gamma,R)$.

The  line complexes $\L_D$ are Schubert varieties \cite{Ful1,Kleim}. Generally the various $\L_D$ lie in special position, so the precise structure of ${V}_\Gamma(A)$ is not accessible via Schubert calculus. Nevertheless,  the dimension and degree of $\pV(\Gamma,R)$ are given in most cases by simple applications of the Pieri rule. We develop these formulae at the end of Section~\ref{sect-lines}.

In Section~\ref{sect-examples} study exhibit some interesting examples from this point of view. We find that the Hessian arrangement, for $R=\overline{\Z}_3,$  has a resonance component which is an irreducible cubic hypersurface in $\P^4$. The arrangement $\overline{\D}$ of directrices consists of four planes in interesting special position in $\P^4$. See Example \ref{hessian}. This may give rise, via exponentiation, to a counter-example to Arapura's result \cite{Arap} for fields of positive characteristic.

We also show that the deleted $B_3$ arrangement of \cite{Suc} has resonance components over fields of characteristic two which do not intersect trivially. Their intersection corresponds to a special point on the translated 1-torus in the first characteristic variety. See Example \ref{deleted1}.

In case $R$ is not an integral domain, the relation determined by a pair of resonant weights on their support may fail to be transitive.
In Example~\ref{db3-z4} we show that the deleted $B_3$ arrangement has such pairs of resonant weights over $\Z_4$, supported on the whole arrangement. Indeed this arrangement has no neighborly partition. In addition, such weights need not have vanishing local weight sums, a necessary condition when $R$ is a domain. On the other hand, these pairs of weights  exponentiate to points which lie on (and thus determine) the translated 1-torus in the characteristic variety. Our treatment of the most general case in the next section is mainly for the purpose of understanding this example.

\end{section}

\begin{section}{Resonance varieties over commutative rings}
\label{resonance}

Let $\A=\{H_1,\ldots,H_n\}$ be an arrangement of distinct linear hyperplanes in $\C^\ell$. The combinatorial structure of \A\ is recorded in the underlying matroid $\sG=\sG(\A).$ This is the matroid on $[n]:=\{1,\ldots, n\}$  whose set of circuits $\Cir$ consists of the minimal sets $C\subseteq [n]$ satisfying $\codim \bigcap_{i\in C} H_i<|C|$. Note that $\sG$ has no circuits of size one or two. From the topological standpoint, we are mainly interested in the topology of the complement $X(\A)=\C^\ell - \bigcup_{i=1}^n H_i$, which is determined to a large, albeit ultimately unknown extent by the underlying matroid $\sG$.

Let $R$  be a commutative ring with $1$. Let $E_R(n)$ denote the free graded exterior algebra over $R$ generated by $1$ and degree-one elements $e_i$ for $i\in [n]$. 

\begin{definition}
The {\em Orlik-Solomon algebra} $A_R(\sG)$ of $\sG$ is the quotient of $E_R(n)$ by the homogeneous ideal $$I=(\partial e_C \ | \ C\in \Cir),$$ where $\partial$ is the usual boundary operator: $\partial e_C=\sum_{k=1}^p (-1)^{k-1}e_{i_1}\we\cdots\hat{e}_{i_k}\cdots \we e_{i_p},$ for $C=\{i_1,\ldots,i_p\}$. 
\end{definition}
The image of $e_i$ in $A_R(\sG)$ is denoted $a_i$. Then $A_R(\sG)$ is a graded-commutative $R$-algebra, generated by $1$ and the degree-one elements $a_i, 1\leq i\leq n.$ According to \cite{OT}, $A^d_R(\sG)$ is a free $R$-module whose rank is independent of $R$. More precisely, the rank of $A^d_R(\sG)$ is equal to the $d^{\rm th}$ Whitney number of the lattice of flats of the matroid $\sG$. The Orlik-Solomon algebra is isomorphic to the cohomology ring of the complement $X(\A)$ with coefficients in $R$,  in case $R$ is an integral domain \cite{OS1,OT}. The generators $a_i$ correspond to logarithmic 1-forms $d\phi_i/\phi_i$ where $\phi_i: \C^\ell \to \C$ is a linear form with kernel $H_i$. 

Resonance varieties over fields were introduced in \cite{F7}; alternate definitions are given in \cite{CS3} and elsewhere. Generalizing the notion to arbitrary commutative rings turns out to be a somewhat treacherous task. We intersperse the definitions given below with some remarks and observations meant to illustrate the perils and justify our conventions.

\begin{definition} Two vectors   $\bx, \bn \in R^n$ are {\em parallel} if the $2 \times 2$ minor determinants $\begin{vmatrix}\x_i&\n_i\\ \x_j&\n_j\end{vmatrix}$ of $\begin{bmatrix}\bx|\bn\end{bmatrix}$ vanish, for $1\leq i<j\leq n.$ 
\end{definition}

If $R$ is a field then \bx\ and \bn\ are parallel precisely when one is a multiple of the other.
The  general case is more subtle, as illustrated by the following theorem. We give the elementary proof; the analogous result for arbitrary linear systems is given in \cite[Section 4.2]{AW}.

\begin{theorem} Let $\bx,\bn \in R^n.$ 
\begin{enumerate}
\item \bx\ and \bn\ are linearly dependent over $R$ if and only if the $2\times 2$ minors of $\begin{bmatrix}\bx|\bn\end{bmatrix}$ have a common nonzero annihilator.
\item If $a\bx+b\bn=0,$ then $a$ and $b$ annihilate all the $2\times 2$ minors of $\begin{bmatrix}\bx|\bn\end{bmatrix}$.
\item If \bx\ and \bn\ are parallel, then \bx\ and \bn\ are linearly dependent over $R$.
\item If $R$ is an integral domain, and \bx\ and \bn\ are linearly dependent over $R$, then \bx\ and \bn\ are parallel.
\end{enumerate}
\label{ann}
\end{theorem}

\begin{proof} The second assertion can be proved by simple Gaussian elimination. This proves necessity in (i). For the converse, note that, if \bx\ and \bn\ are killed by a nonzero element of $R$, then the conclusion holds. Otherwise, we may assume there exists $r\in R-\{0\}$ such that $r$ annihilates the $2\times 2$ minors of $\begin{bmatrix}\bx|\bn\end{bmatrix}$, and $(r\n_1,r\x_1)\not = (0,0)$. Then $(-r\n_1)\bx+(r\x_1)\bn=0$ is a nontrivial dependence relation. Statements (iii) and (iv) are easy consequences.
\end{proof}

For $\bx\in R^n$ set $a_\bx=\sum_{i=1}^n \x_i a_i$. Let $$Z(\bl,R)=\{\be \in R^n \ | \ a_\bl \we a_\be =0\}.$$ Usually $Z(\bl,R)$ is abbreviated to  $Z(\bl)$, when no ambiguity results. For general $R$,  $Z(\bl)$ may not be  a free or even finitely generated module. If $\be\in Z(\bl)$ is not parallel to \bl, we call $(\bl,\be)$ a {\em resonant pair}. The {\em support} $\supp(\bl,\be)$ of a resonant pair  is the set $\{i \in [n] \ | \ \l_i\not = 0 \ \text{or} \ \e_i\not = 0\}.$

\begin{definition} The (degree-one) resonance variety of $A$ is $$\R^1(A)=\{\bl\in R^n \ | \ \exists \, \be \ \text{such that} \ (\bl,\be) \ \text{is a resonant pair}\}.$$
\label{goodres}
\end{definition}

Let $d_\bl: A \to A$ be defined by $d_\bl(x)=a_\bl \we x$. Then $d_\bl^2=0$. Let $H^*(A,d_\bl)$ denote the cohomology of $A$ relative to  $d_\bl$. If $R$ is a field, then $Z(\bl)$ and $H^*(A,d_\bl)$ are $R$-vector spaces and $\dim_R H^1(A,d_\bl)\geq d$ if and only if $\dim_R Z(\bl)\geq d+1.$ In particular, \bl\ is resonant if and only if $H^1(A,d_\bl)\not = 0.$ The last statement is trivially false if $R$ is not a field. Indeed, if $R=\Z$ and \be\ is a nonzero vector, then $\be$ represents a nontrivial class in $H^1(A,2\be)$, but $(2\be,\be)$ is not a resonant pair. 

\begin{definition} Suppose $R$ is an integral domain, with field of quotients $F$. Let $$\R^1_k(A)=\{\bl \in R^n \ | \ \dim_F Z(\bl,F)\geq k+1\}.$$
\label{higher}
\end{definition}

The sets $\R^1(A)$ and $\R_k^1(A)$ are ultimately determined by the arrangement \A\ or matroid $\sG$ and the ring $R$; we will often emphasize this dependence by writing, e.g., $\R^1(A)=\R^1(\sG,R)$.

\begin{theorem} Suppose $R$ is an integral domain. Then $$\R_1^1(\sG,R)=\R^1(\sG,R).$$
\label{gooddef}
\end{theorem}

\begin{proof} Let $F$ be the field of quotients of $R$. Suppose $\dim_F Z(\bl,F)\geq 2.$ Then there is a vector $\be'\in F^n$ such that $(\bl,\be')$ is an $F$-resonant pair. Then there is a multiple $\be$ of $\be'$ that lies in $R^n,$ such that $(\bl,\be)$ is an $R$-resonant pair. Thus $\bl\in \R^1(\sG,R)$. Conversely, if $(\bl,\be)$ is an $R$-resonant pair, then $\{\bl,\be\}$ is a linearly independent subset of $Z(\bl,F)$, so $\bl\in \R_1^1(A,R).$
\end{proof}

\begin{remark}There does not seem to be an extension of Definition \ref{higher} to arbitrary commutative rings for which Theorem \ref{gooddef} remains valid.
If $R$ is a domain, then $\bl \in \R^1_k(A,R)$ if and only if all $(n-k) \times (n-k)$ minors of $d_\bl: A^1 \to A^2$ vanish, i.e., the determinantal rank of $d_\bl$ is at most $n-k-1.$ So in this case the resonance variety $\R_k^1(A,R)$ is indeed an algebraic variety (or affine scheme), defined by the $k^{\rm th}$ Fitting ideal of $d_\bx$ in $R[\x_1,\ldots, \x_n]$. For general $R$,  we can make no analogous statement.
Theorem 2.25 of \cite{AW} (which we specialized in Theorem \ref{ann}) cannot be sharpened in any satisfying way: one can produce $(d+1)$ distinct nonzero elements of $Z(\bl)=\ker(d_\bl)$ provided the $(n-d) \times (n-d)$ minors of $d_\bl$ have a common nonzero annihilator. But one cannot conclude that the resulting set is linearly independent, or even contains two non-parallel vectors, nor that \bl\ is in its span. In particular, the vanishing of all $(n-1)\times(n-1)$ minors of $d_\bl$ does not guarantee that \bl\ is  resonant according to  Definition \ref{gooddef}.

This complication with Fitting varieties and rank of modules over non-domains is obscured in the statement of Theorem 4.5 of \cite{CO}. According to the proof, the quantity $\rank_{\Z_N} H^q(A_{\Z_N}, d_\bl)$ appearing in the statement of that theorem should be interpreted solely in terms of the determinantal rank of $d_\bl: A^q\to A^{q+1}$ when $N$ is not prime.
\label{norank}
\end{remark}
\end{section}

\begin{section}{Combinatorial decomposition of $\R^1(A)$} 
\label{decomp}
In this section we establish algebraic conditions for a pair $(\bl,\be)$ to be resonant, for arbitrary $R$. We define the graph associated with a resonant pair, and thus obtain a decomposition of  the resonance variety $\R^1(A).$ The graph of a resonant pair is shown to be neighborly.  If $R$ is an integral domain, we obtain a similar decomposition of $\R_k^1(\A)$ for all $k\geq 1$. In this case neighborly graphs can be replaced by neighborly partitions of submatroids of \sG. 
This section amounts to a refinement and generalization of the main algebraic results of \cite{F7}. 

Let \A\ be an arrangement, with matroid \sG\ and  set of circuits $\Cir,$ as before. A subset $S\subseteq [n]$ is {\em closed} in $\sG$ if and only if $C-\{i\}\subseteq S$ implies $C\subseteq S$, for all $C\in \Cir$ and $i\in C$. The {\em closure} of $S$ is the (well-defined) smallest closed set containing $S$. The {\em rank} of \sG\ is the size of a minimal set with closure equal to $[n]$.  This is equal to the codimension of $\bigcap_{i=1}^n H_i.$  A {\em line} in $\sG$ is the  closure of a two-point subset of $[n].$ Thus a line in $\sG$ corresponds to a maximal subarrangement of \A\ intersecting in a codimension-two subspace.
A line $X$ in $\sG$ is {\em trivial} if $|X|=2$. We denote the set of lines in $\sG$ by $\X(\sG)$, and the set of nontrivial lines by $\X_0(\sG)$. 

For $\bx=(\x_1,\ldots,\x_n) \in R^n$ and $S\subseteq [n]$, we define the {\em restriction} $\bx_S$ of \bx\ to $S$ to be the element $(\x_i \ | \ i\in S)$ of $R^{|S|}$. The coefficient sum $\sum_{i\in S} \x_i\in R$ is written $\x_S$. 

\subsection*{The rank-two case}
Our approach is based on a characterization of resonant weights in rank two. If $R$ is a domain, and $\bl\in \R^1(\sG,R)$, then $\l_{[n]}=0$, as originally shown in \cite{Y2}. This is false in general, even for \sG\ of rank two. If \sG\ has rank two, the converse holds; the proof of this fact in \cite{F7} carries through for arbitrary commutative rings. A refinement of that argument yields the following result.

\begin{theorem} Suppose $\sG$ has rank two, and $\bl,\be \in R^n$. Then the  following are equivalent:
\begin{enumerate}
\item $a_\bl\we a_\be=0$ 
\item $\begin{vmatrix}\l_{[n]}&\e_{[n]}\\ \l_k&\e_k\end{vmatrix}=0$ for $1\leq k\leq n$.
\item $\l_{[n]}\be=\e_{[n]}\bl$
\end{enumerate}
If {\rm (i)-(iii)} hold, then $\l_{[n]}$ and $\n_{[n]}$ each annihilate the $2\times 2$ minors of $\begin{bmatrix}\bl | \be\end{bmatrix}$.
\label{ranktwo}
\end{theorem}

\begin{proof} First we prove (i) and (ii) are equivalent. The argument relies on two elementary facts  \cite{OT}: $\{a_1\we a_k \ | \ 2\leq k\leq n\}$ forms a basis for $A^2,$ and $(a_i-a_1)\we(a_j-a_1)=0$ for $1\leq i, j\leq n$. Then 
\begin{equation*}
\begin{split} a_\bl\we a_\be&=\biggl(\sum_{i=1}^n \l_i a_i\biggr)\we\biggl(\sum_{j=1}^n \e_j a_j\biggr)\\
&=\biggl(\l_{[n]}a_1+\sum_{i=2}^n \l_i (a_i-a_1)\biggr)\we\biggl(\e_{[n]}a_1+ \sum_{j=2}^n \e_j(a_j-a_1)\biggr)\\
&= \sum_{k=2}^n \begin{vmatrix}\l_{[n]}&\e_{[n]}\\ \l_k&\e_k\end{vmatrix} a_1\we a_k,
\end{split}
\end{equation*}
where we have used the second elementary fact. Then $a_\bl\we a_\be=0$ if and only if $\begin{vmatrix}\l_{[n]}&\e_{[n]}\\ \l_k&\e_k\end{vmatrix}=0$ for $2\leq k\leq n$ by the first elementary fact. This statement is equivalent to (ii).

Since (ii) says $\l_{[n]}\n_k=\n_{[n]}\l_k$ for every $k$, (ii) is equivalent to (iii). The final assertion follows from Theorem \ref{ann}(ii).
\end{proof}

When $R$ is a domain we obtain the usual characterization of rank-two resonant weights.
\begin{corollary} If $R$ is a domain, then $a_\bl\we a_\be=0$ if and only if \bl\ and \be\ are parallel, or $n\geq 3$ and $\bl_{[n]}=0= \be_{[n]}$.
\label{usual}
\end{corollary}

There is no analogue of Corollary \ref{usual} for arbitrary rings. The coefficient sums $\l_{[n]}$ and $\e_{[n]}$ for a resonant pair $(\bl,\be)$ need not vanish in general: for example, take $\bl=(-1,3,1)$ and $\be=(-1,1,3)$ over $R=\Z_6.$ If $(\bl,\be)$ is an $R$-resonant pair, then $n\geq 3$ and $\l_{[n]}$ and $\e_{[n]}$, if not zero, must be zero-divisors, by Theorem \ref{ranktwo} and Theorem \ref{ann}. In this case \bl\ and \be\ are linearly dependent over $R$. The converse does not hold. The next results detail the conclusions one may draw in the general case.

Henceforth we will abuse the standard terminology by saying $s\in R$ is a {\em zero divisor} if there exists $r\in R-\{0\}$ such that $rs=0,$ that is, if $\Ann_R(s)\not =0.$ Recall $Z(\bl)=\{\be\in R^n \ | \ a_\bl\we a_\be=0\}$. Let $\Delta=\{\bx\in R^n \ | \ \x_{[n]}=0]\}.$ 

\begin{corollary}  $\Ann_R(\bl)^n \cup (\Ann_R(\l_{[n]})^n\cap \Delta)\subseteq Z(\bl).$
\end{corollary}

Every vector in $\Ann_R(\bl)^n$ is parallel to \bl. Also $\Ann_R(\bl)$ is contained in $\Ann_R(\l_{[n]})$; if the inclusion is proper we can construct a resonant pair for \bl, using a variant of the proof of Theorem \ref{ann}. This is the closest we can get to a generalization of Corollary \ref{usual}.

\begin{theorem}  Suppose $\sG$ has rank two, and $n\geq 3$. Let $\bl \in R^n$.
If $\l_{[n]}$ is  a zero divisor, with $\Ann_R(\l_{[n]})\not = \Ann(\bl),$ then \bl\ is resonant.
\label{closest} 
\end{theorem} 

\begin{proof}  Let  $r\in \Ann_R(\l_{[n]})-\Ann_R(\bl).$ Then  $r\l_{[n]}=0$ and $r\l_i\not=0$ for some $i$; without loss $i=1.$ Since $n\geq 3,$ we may set $\be=(0,r,-r,0,\ldots,0).$ Then $\l_{[n]}\be=\e_{[n]}\bl$ and  $\begin{vmatrix}\l_1&\e_1\\ \l_2&\e_2\end{vmatrix}\not =0.$ Thus $(\bl,\be)$ is a resonant pair.
\end{proof}

The relation $\l_{[n]}\be=\e_{[n]}\bl$  can be interpreted in terms of linear line complexes as in Section \ref{sect-lines}. See Remark \ref{ranktwocomplex}. This raises the possibility that $\R^1(\sG,R)$ may be nonlinear if $R$ is not a domain, even for \sG\ of rank two.

\subsection*{The general case}

Using the grading of $A^2$ by $\X(\sG)$, Theorem \ref{ranktwo} yields a characterization of resonant weights for matroids of any rank.
\begin{theorem} Suppose $\sG$ is a matroid of arbitrary rank. Then $\be\in Z(\bl)$ if and only if, for every $X\in \X(\sG)$, either \begin{enumerate}
\item $\bl_X$ and $\be_X$ are parallel, or
\item $X\in \X_0(\sG)$ and  $\l_X\be_X=\e_X\bl_X.$ 
\end{enumerate}
\label{res-pair}
\end{theorem}

\begin{proof} Note that $a_\bl\we a_\be=\sum_{i<j} \begin{vmatrix}\l_i&\e_i\\ \l_j&\e_j\end{vmatrix}  \ a_i\we a_j.$ There is a direct sum decomposition $$A^2=\underset{X\in \X(\sG)}{\oplus} A^2_X,$$ where $A_X$ is the subalgebra of $A$ generated by $\{a_i \ | \ i\in X\}$ \cite{OT}. Since $X\in \X(\sG)$ is a rank-two submatroid of $\sG$, Theorem \ref{ranktwo} yields the result.  
 \end{proof}

In case (ii) above $\l_X$ and $\e_X$ are zero divisors.
We make no claim about $\l_{[n]}$ in general. Example \ref{db3-z4} exhibits a resonant weight \bl\ over $\Z_4$, supported on a matroid of rank 3, for which $\l_X\not = 0$ for some $X\in \X_0(\sG).$ In this example $\l_{[n]}=0.$  We do not know whether $\l_{[n]}$ must be a zero divisor for resonant \bl\ if \sG\ has rank greater than two, though it cannot be a unit by \cite{Y2}.

\subsection*{Neighborly graphs and partitions}
Suppose $(\bl,\be)$ is a resonant pair. Define a graph $\Gamma=\Gamma_{(\bl,\be)}$ with vertex set  $[n]$, and $\{i,j\}$ an edge of $\Gamma $ if and only if  $\begin{vmatrix}\l_i&\e_i\\ \l_j&\e_j\end{vmatrix}=0$. If $\{i,j\}$ is an edge of $\Gamma$ we write $\{i,j\}\in\Gamma$. Note that $\{i,j\}\in \Gamma$ for every trivial line $\{i,j\}\in \X(\sG)$. Also, if $i\not \in \supp(\bl,\be)$ then $i$ is a cone vertex in $\Gamma$, adjacent to every other vertex. A {\em clique} in $\Gamma$ is a set of vertices contained in a complete subgraph. Since \bl\ and \be\ are not parallel, $[n]$ itself is not a clique of $\Gamma$. By Theorem \ref{res-pair} and the definition of $\Gamma_{(\bl,\be)}$, we have the following.
 
\begin{corollary} If $X\in \X_0(\sG)$ is not  a clique of $\Gamma_{(\bl,\be)}$, then $\l_X\be_X=\e_X\bl_X$.
If $X\in \X(\sG)$ is a clique in $\Gamma_{(\bl,\be)}$, then $\bl_X$ is parallel to $\be_X$.
\end{corollary}

We define a {\em block} of $\Gamma$ to be a maximal clique. The blocks of $\Gamma$ cover $[n]$, but need not be disjoint. A cone vertex of $\Gamma$ is contained in every block.
\begin{definition} A graph $\Gamma$ with vertex set $[n]$ is {\em neighborly,} or more precisely {\em $\sG$-neighborly,} if for every $X\in \X(\sG)$ and every block $S$ of $\Gamma$, $|X\cap S|\geq |X|-1$ implies $X\subseteq S$.
\label{neighborlydef}
\end{definition}

Observe that a \sG-neighborly graph must include among its edges all the trivial lines $\{i,j\}\in \X(\sG)$. Also, if $i$ is a cone vertex of $\Gamma$, then $\Gamma$ is $\sG$-neighborly if and only if the induced subgraph on $[n]-\{i\}$ is $(\sG-i)$-neighborly.

\begin{theorem} Let $\Gamma=\Gamma_{(\bl,\be)}$ be the graph associated with a resonant pair of weights. Then $\Gamma$ is neighborly.
\label{clique}
\end{theorem}

\begin{proof}
Let $X\in \X(\sG),$ and $i\in X$ with $X-\{i\}$ a clique of $\Gamma.$  Suppose $X$ is not  a clique. Then $\l_X\e_k=\e_X\l_k$ for $k \in X$ by Theorem \ref{res-pair}(ii). But $\l_j\e_k=\e_j\l_k$ for $j,k\in X-\{i\}$ since $X-\{i\}$ is a clique. We conclude $\l_i\e_k=\e_i\l_k$ for all $k\in X$, a contradiction. 
\end{proof}

We say the graph $\Gamma$ is {\em transitive} if $\{i,j\}, \{j,k\}\in \Gamma$ implies $\{i,k\}\in \Gamma$. If $\Gamma$ is a transitive graph then the components of $\Gamma$ are cliques, hence are the blocks of $\Gamma$. A transitive graph with a cone vertex is a complete graph. If $\Gamma$ is a transitive neighborly graph with no cone vertices, the blocks of $\Gamma$ form a {\em neighborly partition} of $\sG$ in the sense of \cite{F7}. 

\begin{theorem} Suppose $R$ is an integral domain, and $\Gamma=\Gamma_{(\bl,\be)}$ is the graph of an $R$-resonant pair. Then 
\begin{enumerate}
\item $\supp(\bl,\be)$ coincides with the set of non-cone vertices of $\Gamma,$ and
\item $\Gamma$ is transitive on $\supp(\bl,\be).$
\end{enumerate}
\label{eqrln}
\end{theorem}
\begin{proof} Suppose $\{i,j\}, \{j,k\}\in \Gamma$, with $j\in \supp(\bl,\be)$. Then $\l_i\e_j=\l_j\e_i$ and $\l_j\e_k=\l_k\e_j$.  Since $j\in \supp(\bl,\be)$ we can assume without loss that $\e_j\not =0$.
We have $\l_i\e_j\e_k=\l_j\e_i\e_k=\e_i\l_k\e_j$, which then implies $\l_i\e_k=\l_k\e_i$ since $R$ is a domain. Thus $\{i,k\}\in \Gamma$. This proves the second assertion, and also shows that a cone vertex cannot lie in $\supp(\bl,\be)$, else $\Gamma$ itself is a clique. 
\end{proof}
In Example \ref{db3-z4} we will see  a resonant pair $(\bl,\be)$ over $R=\Z_4$ for which $\supp(\bl,\be)$ includes some cone vertices of $\Gamma_{(\bl,\be)}$. In particular, $\Gamma_{(\bl,\be)}$ is not transitive on $\supp(\bl,\be)$.

\subsection*{Combinatorial components}
For $\Gamma$ an arbitrary graph on $[n]$, set $$\X_\Gamma(\sG)=\{X\in \X_0(\sG) \ | \ X \ \text{is not a clique of} \ \Gamma\}$$ and 
$$K(\Gamma,R)=\{\bx \in R^n \ | \ \x_X  \ \text{is a zero divisor for every} \ X\in \X_\Gamma(\sG)\}.$$ 

For $\bl\in K(\Gamma,R)$ we define  $Z_\Gamma(\bl,R)$ to be the set of those $\be\in K(\Gamma,R)$ satisfying
\begin{enumerate}
\item $\l_X\be_X=\e_X\bl_X$ for all $X\in \X_\Gamma(\sG)$, and
\item$\begin{vmatrix}\l_i&\l_j\\ \e_i&\e_j\end{vmatrix} =0$ for every edge $\{i,j\}$  of $\Gamma$.
\end{enumerate}
In particular, if $\be\in Z_\Gamma(\bl,R)$, and $S$ is a clique of $\Gamma$, then $\be_S$ is parallel to $\bl_S$. The converse may not be true, that is, the graph $\Gamma_{(\bl,\be)}$ may have more edges than the original graph $\Gamma$. We will write $Z_\Gamma(\bl,R)$ as $Z_\Gamma(\bl)$ when it is not ambiguous.
If $R$ is a domain, then $$K(\Gamma,R)=\{\bx \in R^n \ | \ \x_X=0 \ \text{for every} \ X\in \X_\Gamma(\sG)\},$$  and condition (i) is vacuous. In this case $Z(\bl,R)$ is a submodule of $R^n$.

\begin{corollary} $Z_\Gamma(\bl)\subseteq Z(\bl)$.
\label{graph-cycles}
\end{corollary}

\begin{proof}
Let $\be\in Z_\Gamma(\bl).$ If $\bl_X$ and $\be_X$ are not parallel, then $X\in \X_\Gamma(\sG)$ and $\l_X\be_X=\e_X\bl_X$. Then $\be\in Z(\bl)$ by Theorem \ref{res-pair}.
\end{proof}

\begin{definition} The {\em combinatorial component} of $\R^1(\sG,R)$ corresponding to a graph $\Gamma$ is $$V^1(\Gamma,R)=\{\bl \in K(\Gamma,R) \ | \  \exists \, \be \in Z_\Gamma(\bl) \ \text{such that \be\ is not parallel to} \  \bl\}.$$
\end{definition}
In case $R$ is an integral domain, $V^1(\Gamma,R)=\{\bl \in K(\Gamma,R) \ | \ \dim_F Z(\bl,F)\geq 2\},$ where $F$ is the quotient field of $R$.
Let $\NG(\sG)$ denote the set of $\sG$-neighborly graphs with vertex set $[n]$. 
Let $\NG(\sG,R)$ denote the set of $\Gamma\in \NG(\sG)$ for which $K(\Gamma,R)$ contains a pair of non-parallel vectors. If $R$ is a domain, then $\NG(\sG,R)=\{\Gamma \in \NG(\sG) \ | \ \dim_F K(\Gamma,F)\geq 2\},$ with $F$ as above.

We can now establish the decomposition theorem for general $R$.

\begin{theorem} For any commutative ring $R$, $$\R^1(\A,R)=\bigcup_{\Gamma\in \NG(\sG,R)} V^1(\Gamma,R).$$
\label{bigthm3}
\end{theorem}
\begin{proof} Suppose $\bl\in \R^1(\A,R)$. Then there exists $\be\in R^n$ such that $(\bl,\be)$ is a resonant pair. Let $\Gamma=\Gamma_{(\bl,\be)}.$ Then $\Gamma\in \NG(\sG)$ by Theorem \ref{clique}. Furthermore, $Z_\Gamma(\bl)$  contains the non-parallel elements \bl\ and \be. Thus $\Gamma\in \NG(\sG,R)$ and $\bl\in V^1(\Gamma,R)$. The other inclusion holds by Corollary \ref{graph-cycles}.
\end{proof}

Suppose $R$ is an integral domain. Let $\NP(\sG)\subseteq \NG(\sG)$  be the set of $\sG$-neighborly graphs which are transitive on non-cone vertices. Thus $\NP(\sG)$ corresponds to the set of neighborly partitions of submatroids of $\sG$. Let $\NP(\sG,R)=\NP(\sG)\cap \NG(\sG,R).$  As noted above,  $K(\Gamma,R)$ is the kernel of a row-selected submatrix of the $|\X_0(\sG)| \times n$ point-line incidence matrix of $\sG$;  $\Gamma\in \NP(\sG)$ lies in $\NP(\sG,R)$ if and only if this matrix has nullity at least two.

\begin{theorem} For any integral domain $R$, $$\R^1(\A,R)=\bigcup_{\Gamma\in \NP(\sG,R)} V^1(\Gamma,R).$$
\end{theorem}
\begin{proof} In the proof of Theorem \ref{bigthm3}, the graph $\Gamma=\Gamma_{(\bl,\be)}$ is transitive on the non-cone vertices of $\Gamma$ by Theorem \ref{eqrln}. 
\end{proof}

The term ``component"  is potentially misleading: $V^1(\Gamma,R)$ may be trivial for some graphs $\Gamma\in \NG(\sG,R)$. Furthermore, we make no claim that $V^1(\Gamma,R)$ is irreducible, even for $R$ an algebraically closed field, although we have no examples to the contrary.  Without more precise information about the incidence structure of $\sG$, the most one can say in this regard is the following.

\begin{theorem} Suppose $\Gamma'$ is a subgraph of $\Gamma$ and $\X_{\Gamma'}(\sG)\subseteq \X_\Gamma(\sG)$. Then $V^1(\Gamma,R)\subseteq V^1(\Gamma',R)$.
\label{nested}
\end{theorem}

Observe that the two conditions in Theorem \ref{nested} are somewhat in opposition: the fewer edges in $\Gamma$, the more lines (potentially) in $\X_\Gamma(\sG)$. This tension, together with the neighborliness required of $\Gamma$, accounts for the dearth of matroids supporting resonant pairs over integral domains.

The {\em support} of $V^1(\Gamma,R)$  is the set of indices $i$ such that $\l_i\not =0$ for some $\bl\in V^1(\Gamma,R)$. In case $R$ is an integral domain, the support of $V^1(\Gamma,R)$ is the set of non-cone vertices of $\Gamma$, by Theorem \ref{eqrln}. We say $V^1(\Gamma,R)$ is {\em essential} if its support is $[n]$.

\subsection*{Higher order resonance varieties}
Now suppose $R$ is an integral domain, so that $\R^1_k(\sG,R)$ is defined. For $\bl\in R^n$ we define a single graph $\Gamma$, depending only on $\bl$, such that $\be\in Z_\Gamma(\bl)$ for every resonant pair $(\bl,\be)$. 

Let $E$ be a field extension of $R$. An element $\bm \in E^n$ is called a {\em generic partner} of \bl\ if the following conditions are satisfied:

\begin{enumerate}
\item if $X\in \X_0(\sG)$ and $\m_X=0$, then $\e_X=0$ for all $\be \in Z(\bl,R)$, and
\item if $1\leq i<j\leq n$ and $\begin{vmatrix}\l_i&\l_j\\ \m_i & \m_j\end{vmatrix}=0$, then $\begin{vmatrix}\l_i&\l_j\\ \e_i & \e_j\end{vmatrix}=0$ for every $\be\in Z(\bl,R)$.
\end{enumerate}

Every $\bl\in R^n$ has a generic partner: let $E$ be the algebraic closure of the quotient field of $R$, and apply Hilbert's Nullstellensatz. A more precise existence theorem may be useful for computational purposes. 
\begin{theorem}If $E$ is a field extension of $R$ satisfying
$$|E|>\binom{n}{2}+|\X_0|,$$ then every $\bl\in R^n$ has a generic partner in $E^n$. If $R$ is infinite,  \bl\ has a generic partner in $R^n.$
\end{theorem}

\begin{proof}
Let $d=\dim_E Z(\bl,E).$ The linear equations $\x_X=0$ determine $|\X_0|$ hyperplanes in $Z(\bl,E)$. For $\bl$ fixed, the equations $\begin{vmatrix}\l_i&\l_j\\ \x_i & \x_j\end{vmatrix}=0$ also define hyperplanes in $Z(\bl,E)$, at most $n \choose 2$ of them. Since
$$|Z(\bl,E)|=|E|^d>|E|^{d-1}\bigl({n\choose 2}+|\X_0|\bigr),$$
there is a point $\bm$ of $Z(\bl,E)$ missing the aforementioned subspaces.

If $R$ is infinite take $E$ to be the quotient field of $R$. Then \bl\ has a generic partner $\bm\in E^n,$ by preceding argument. Some nonzero multiple of \bm\ will lie in $R^n,$ and remains a generic partner of \bl.
\end{proof}

We define $\Gamma_\bl=\Gamma_{(\bl,\bm)}$, where \bm\ is a generic partner of \bl. It follows from condition (ii) that $\Gamma_\bl$ is well-defined. By Theorem \ref{eqrln}, $\Gamma_\bl \in \NP(\sG,E)$.
Note $\supp(\bl,\bm)\supseteq\supp(\bl,\be)$ for every $\be\in Z(\bl,R)$.

\begin{theorem} Suppose $(\bl,\be)$ is a resonant pair over a domain $R$. Then $\Gamma_\bl\in \NP(\sG,R)$ and $\be \in Z_{\Gamma_\bl}(\bl,R)$. 
\label{comb-cycles}
\end{theorem}

\begin{proof} If $X\in \X_{\Gamma_\bl}(\sG),$ then $\bl_X$ and $\bm_X$ are not parallel, so $\l_X=0=\m_X$ by Theorem \ref{res-pair}. 
Then $\e_X=0$ by genericity of \bm. Hence $\be\in K(\Gamma_\bl,R)$.
Also, by (ii) above,  $\begin{vmatrix}\l_i&\e_i\\ \l_j&\e_j\end{vmatrix}=0$ for every $\{i,j\}\in \Gamma_\bl$. Thus $\be\in Z_{\Gamma_\bl}(\bl)$. Since $\bl, \be \in K(\Gamma,R)$ are not parallel, $\Gamma_\bl\in \NP(\sG,R).$
\end{proof}
By Theorem \ref{comb-cycles} and Corollary \ref{graph-cycles}, we have 
\begin{corollary}  $Z(\bl,R)=Z_{\Gamma}(\bl,R)$  for $\Gamma=\Gamma_\bl$.
\label{goodgraph}
\end{corollary}
Let $V^1_k(\Gamma,R)=\{\bl \in R^n \ | \ \dim_F Z_\Gamma(\bl)>k\}.$ Then  $V_k^1(\Gamma,R)=0$ unless $\dim_F K(\Gamma,F)>k$. Note that $V_1^1(\Gamma,R)=V^1(\Gamma,R)$.
\begin{corollary} Suppose $R$ is an integral domain.
Then $$\R^1_k(\sG,R)=\bigcup_{\Gamma\in \NP(\sG,R)} V^1_k(\Gamma,R).$$
\label{bigthm2}
\end{corollary}

\begin{proof} By Corollary \ref{goodgraph}, if $\bl\in \R_d(\sG,R)$, then $\bl\in V_d(\Gamma_\bl,R)$. The reverse inclusion follows from Corollary \ref{graph-cycles}.
\end{proof}

\end{section}

\begin{section}{The structure of $ V^1(\Gamma,R)$}
\label{sect-lines}
Geometers of the early $20^{\rm th}$ century understood well the connection between skew-symmetric forms (over $\Re$ or $\C$), null polarities on projective space, and projective line complexes \cite{Veblen-Young,Pottman}.  In this section we return to their methods in our more general setting.

Throughout this section, we assume for simplicity that $R$ is an algebraically closed field. Since $V^1(\Gamma,R)$ is preserved under the diagonal action of $R^*$ we consider its projective image.
We will see that this projective variety is the carrier of an algebraic line complex determined by certain projective subspaces associated with $\Gamma$. For background on line complexes, Grassmannians and Schubert varieties we refer the reader to \cite{Pottman,Kleim,Ful1,GH}.

Let $K$ be a vector space of dimension $k>0$ over $R$, and let $\P(K)=\{R\bx \ | \ \bx\in K-\{0\}\}$ be the projective space of $K$. The standard projective space $\P(R^k)$ is denoted $\P^{k-1}$. If $\bx, \bn\in K-\{0\}$ are not parallel, we denote by $\bx\ast\bn$ the line in $\P(K)$ spanned by $R\bx$ and $R\bn$. Thus $\bx\ast\bn=\P(R\bx+R\bn)$. If $\bx\in K$ and  $D$ is a nontrivial subspace of $K$, let $\bx\ast D=\P(R\bx+D);$ if $D$ and $D'$ are subspaces of $K$, let $D\ast D'=\P(D+D')$. We will usually abbreviate $R\bx$ to $\bar{\bx},$ and $\P(D)$ to $\bar{D}$.

\subsection*{Projective line complexes}  A line $L=\bx\ast\bn$ in $\bar{K}$ corresponds to an element of the Grassmanian $\G(2,k)$. The $2 \times 2$ minors $L_{ij}=\x_i\n_j-\x_j\n_i,$ $1\leq i < j \leq d,$ of the matrix $\begin{bmatrix}\bx | \bn\end{bmatrix}$ are called the {\em line coordinates} of $L$. They are determined up to scalar multiple by $L$, independent of the choice of \bx\ and \bn. The Pl\"ucker embedding $L \mapsto [L_{ij} \ : \ 1\leq i<j\leq d]$ identifies $\G(2,k)$ with a $2(k-2)$-dimensional subvariety of $\P^N$, $N=\binom{k}{2} -1.$ The Grassmann-Pl\"ucker relations give a particular set of defining equations for the image $\G(2,k)\subset \P^N$.

A {\em line complex} in $\bar{K}$ is an algebraic subset \L\ of the Grassmannian $\G(2,k)$ under the Pl\"ucker embedding, i.e., a set of lines $\bx \ast\bn$ in $\bar{K}$ given by a system of polynomial equations in the line coordinates $L_{ij}$. The {\em carrier} of a line complex \L\ is the algebraic set $|\L|\subseteq \bar{K}$ of points lying on lines of \L. That is, $|\L|=\bigcup \L$. A {\em variety ruled by lines} is a variety which is the carrier of some line complex.

We are mainly interested in line complexes of the following form. If $D$ is a nontrivial subspace of $K$, set $$\L_D=\{L\in \G(2,k) \ | \ L \cap \bar{D}\not= \emptyset \}.$$ In fact $\L_D$ is a {\em linear line complex}: if $B$ is a matrix whose columns give a basis for $D$, then $L=\bx\ast \bn\in \L_D$ if and only if all maximal minors of  $\begin{bmatrix}B|\bx|\bn\end{bmatrix}$ vanish. Using the Laplace expansion these minors become linear equations in the $L_{ij}$. 
If $\mathcal D$ is an arrangement of nontrivial subspaces in $K$, let $$\L(\D)=\bigcap_{D\in \D} \L_D.$$ 

\subsection*{Combinatorial components as ruled varieties} Let \sG\ be a simple matroid on $[n]$. Fix a graph $\Gamma\in \NP(\sG,R)$ and set $K=K(\Gamma,R)$, as defined in the Section \ref{decomp}. Let $\pV(\Gamma)=\pV(\Gamma,R)$ be the projective image of $V^1(\Gamma,R)$.  Assume $\pV(\Gamma,R)$ is nonempty.

Observe that $\pV(\Gamma)$ is a ruled variety. Indeed, if $\bl\in V^1(\Gamma,R)$, then there exists $\be\in Z_\Gamma(\bl)$ not parallel to \bl. If $\bx\in R\bl+R\be,$ then $Z_\Gamma(\bx)$ contains $R\bl+R\be,$  so $\dim_R  Z_\Gamma(\bx)>1.$ This implies $\bl\ast\be\subseteq \pV(\Gamma).$ 

We proceed to identify the underlying line complex.
Recall that a block of $\Gamma$ is a maximal clique. If $S$ is a block of $\Gamma$, set $$D_S=\{\bx\in K\ | \ \x_i=0 \ \text{for all} \ i\in S\}.$$ The {\em arrangement of directrices} associated with $\Gamma$ is the collection $\D_\Gamma$ of subspaces $D_S,$ where $S$ is a block of $\Gamma$. Note, if $i$ is a cone vertex of $\Gamma$, then every $D\in \D_\Gamma$ is contained in the coordinate hyperplane $\x_i=0$. 

\begin{theorem} Suppose $\bl\in V^1(\Gamma,R)$ and $\be\in Z_\Gamma(\bl)$ is not parallel to \bl. Then, for any block $S$ of $\Gamma$, $\bl\ast\be$ meets $\bar{D}_S$. Conversely, if $L=\bl\ast\be$ is a line in $\bar{K}$ which meets $\bar{D}_S$ for every block $S$ of $\Gamma$, then $L\subseteq \pV(\Gamma)$.
\label{lines}
\end{theorem}

\begin{proof}  Let  $S$ be a block of $\Gamma$. Then $\bl_S$ is parallel to $\be_S$, by definition of $Z_\Gamma(\bl)$. Then we can find scalars $a,b\in R$ such that $\bx=a\bl+b\be$ satisfies $\bx\not=0$ and $\bx_S=a\bl_S+b\be_S=0$.  Also $\bl,\be\in K$ by definition. Then $\bx\in D_S$, so $\bar{\bx}\in (\bl\ast\be) \cap \bar{D}_S$. For the converse, suppose \bl\ and \be\ are not parallel, and $\bl\ast\be$ meets $\bar{D}_S$ for each block $S$ of $\Gamma$. Then for every  block $S$ there exist $a,b\in R$ such that $a\bl+b\be\in D_S$, consequently $a\bl_S+b\be_S=0$. Hence $\be_S$ is parallel to $\bl_S$. Since $\bl,\be\in K$ by assumption, this puts $\be$ in $Z_\Gamma(\bl)$. This implies $\bl\ast\be\subseteq \pV(\Gamma)$ by our previous observation. 
\end{proof}

\begin{corollary} The combinatorial component $\pV(\Gamma)$ is the carrier of the line complex $\L(\D_\Gamma)$.
\label{bigthm1}
\end{corollary}

The argument that  $\pV(\Gamma)$ is a ruled variety and the proof of Theorem \ref{lines} are not valid for rings with zero divisors. 

\subsection*{Schubert calculus in ${\boldsymbol{{\mathcal G}(2,k)}}$} 
The linear line complexes $\L_D$ are in fact Schubert varieties in $\G(2,k)$. The classical intersection theory of Schubert varieties can be used to determine the degree of $|\L(\D_\Gamma)|=\pV(\Gamma)$ in many cases. We remind the reader of
some elementary aspects of the theory as it applies to $\G(2,k)$. See \cite{Kleim,Harris, Ful1,GH} for a more complete development. 

Given a complete flag of subspaces $0=K_0\subset K_1\subset \cdots \subset K_k=K$ and a pair $\sigma=(i_1,i_2)$ of  integers $0\leq i_2\leq i_1\leq k-2$, the associated Schubert variety is the collection $W_\sigma$ of lines $L\in \G(2,k)$ satisfying

\begin{enumerate}
\item $L\cap \P(K_{k-1-i_1})\not = \emptyset,$ and
\item $L\subseteq \P(K_{k-i_2})$.
\end{enumerate}

The pair $(i_1,i_2)$ is usually represented by a Ferrers diagram, or {\em shape}, consisting of a left-justified array of two rows, with $i_1$ boxes in the first row and $i_2$ in the second. The condition on $(i_1,i_2)$ is that this array has non-increasing row-lengths and fits in a $2 \times (k-2)$ rectangle. We call such shapes {\em admissible}.

Two Schubert varieties determined by the same shape, but different flags, are projectively equivalent. The codimension of $W_{(i_1,i_2)}$ in $\G(2,k)$ is $i_1+i_2$. By choosing a flag which includes the subspace $D$, we obtain the following.

\begin{theorem} The line complex $\L_D$ is equivalent to $W_{(s,0)}$ where $s=\codim(D)-1$. The codimension of\/ $\L_D$ in $\G(2,k)$ is $\codim(D)-1$.
\label{cell}
\end{theorem}

The Pieri rule describes the intersection of $W_{(s,0)}$ with $W_\sigma$ as a sum of Schubert varieties, up to rational equivalence. 

\begin{theorem} $$W_{(s,0)}\cdot W_\sigma=\sum W_\tau,$$ where the sum is indexed by those admissible shapes $\tau$ that can be obtained from the shape $\sigma$ by adding $s$ boxes, no two in the same column.
\end{theorem}
The dual Pieri rule in $\G(2,k)$ is somewhat restricted in our special case.
\begin{theorem} $$W_{(1,1)}\cdot W_{(i_1,i_2)}=W_{(i_1+1,i_2+1)},$$ if $(i_1+1,i_2+1)$ is admissible, and vanishes otherwise.
\end{theorem}

These intersection formulae hold in the Chow ring of $\G(2,k)$. They can be interpreted geometrically: if $W$ and $W'$ are Schubert varieties in $\G(2,k)$ which intersect properly, that is, with the expected codimension, then $W\cap W'$ is a union of Schubert varieties  with multiplicities, whose types are given by the terms in the expansion of $W\cdot W'$.  If the intersection is generically transverse, then all multiplicities are equal to one \cite[Section 5.3]{Ful2}. General translates of Schubert varieties in $\G(2,k)$ meet generically transversely, in any characteristic \cite{Sottile}.

\subsection*{Dimension and degree of $\boldsymbol{\pV(\Gamma)}$}
Fix an arrangement $\D$ of nontrivial subspaces of $K$, and let $\L=\L(\D)$. We establish some relations between $\L$ and its carrier $|\L|$. We use \cite{Harris} and \cite{Ful2} as general references  on intersection theory.

If $\bar{\bx} \in |\L|$, the {\em cone} of $\bar{\bx}$ in \L\ is the line complex $\L_\bx=\{L\in \L \ | \ \bar{\bx}\in L\}$.
We define the {\em depth} of $\bar{\bx}$ in \L\ by $\depth(\bx)=\dim |\L_\bx|.$  It is easy to see that $\dim \L_\bx = \depth(\bx)-1.$ If $|L|$ is irreducible, define the depth of \L\ to be the depth of a generic point on $|\L|$. Then $1\leq \depth(\L)\leq \dim \L +1.$

\begin{theorem} If $|L|$ is irreducible, then  $\dim |\L|=\dim \L - \depth(\L) +2.$
\label{dimensions}
\end{theorem}

\begin{proof}  Let ${\mathcal I}$ denote the incidence variety $\{(\bx,L)\in \bar{K}\times \L \ |  \ \bx\in L\}\subseteq \bar{K}\times \G(2,k)$. The fiber of the projection ${\mathcal I} \to \L$ over $L \in \L$ is $L$ itself, of dimension one. Then $\dim {\mathcal I}=\dim \L + 1.$
The fiber of the other projection ${\mathcal I} \to \bar{K}$ over $\bar{\bx} \in \bar{K}$ is $\L_\bx$, and the image of this projection is $|\L|$. Using $\dim \L_\bx=\dim |\L_\bx| -1$,
the result follows.
\end{proof}
In our situation $\L_\bx$ has a special form.
\begin{lemma} $$|\L_\bx|=\bigcap_{D\in \D}\bx\ast D.$$ In particular $|\L_\bx|$ is linear, and $\L_\bx$ consists of the lines through $\bar{\bx}$ in $|\L_\bx|$. 
\label{cone}
\end{lemma}

\begin{proof} If $\bar{\bn}\in |\L_\bx|,$ with $\bar{\bn}\not =\bar{\bx},$ then $\bx\ast\bn$ meets $\bar{D}$ for every $D\in \D,$ so $\bar{\bn}\in \bigcap_{D\in \D}\bx\ast D.$ Conversely, if $\bar{\bn}\in \bx\ast D,$ then $\bar{\bn}$ lies on a line which contains $\bar{\bx}$ and meets $\bar{D}$, equivalently, $(\bx\ast\bn)\cap\bar{D}\not = \emptyset$.
\end{proof}
According to the lemma, $\L_\bx$ is equivalent to the Schubert variety $W_\sigma$ for $\sigma=(k-2,k-1-\depth(\bx)).$ If $\D=\D(\Gamma)$ for some graph $\Gamma$, then, by Theorem \ref{lines}, $|\L_\bl|=\P(Z_\Gamma(\bl)).$  Recall $V^1_k(\Gamma,R)=\{\bl \in V^1(\Gamma,R) \ | \ \dim Z_\Gamma(\bl,R)>k\}.$ We will write $\pV_k(\Gamma)$ for the projective image of $V^1_k(\Gamma,R).$

\begin{corollary} $\bl \in V^1_k(\Gamma,R)$ if and only if $\bar{\bl}$ has depth $k$ in $\L(\D_\Gamma)$.
\end{corollary}

At this point, we have no method to compute any of the dimensions which appear in the Theorem \ref{dimensions} above, except by inspection or computational algebra (e.g., {\em Macaulay2}). The description of $|\L_\bx|$ in Lemma \ref{cone} does not yield a formula for $\depth(\bx)$ because the subspaces $\bx\ast D$ may not be in general position.

For $D$  a subspace of $K$, set $c(D)=\codim(D)-1=\codim(\bar{D})-1$. Then $\L_D$ is projectively equivalent to $W_{(c(D),0)}$. The following observation is adapted from \cite{Pottman}; see also \cite[Example 19.11]{Harris} and \cite[Example 8.3.14]{Ful3}.
\begin{theorem}  Let $\D$ be an arrangement of nontrivial subspaces of $K$, and $\L=\L(\D).$ Suppose $|\L|$ is irreducible. Let $D_0$ be a subspace of $K$ of codimension equal to $\dim |\L|$, with $\bar{D}_0$ in general position relative to $|\L|$. Then the degree of $|\L|$ is given by the following formula: $$W_{(c(D_0),0)}\cdot \L=(\deg |\L|) \,W_\sigma,$$ where $\sigma=(k-2,k-1-\depth(\L)).$
\label{degree}
\end{theorem}

\begin{proof} Let $B= \bar{D}_0\cap|\L|$. The hypotheses on $D_0$ imply that $B$ consists of $\deg |\L|$ points, each of depth equal to $\depth(\L)$. We have $\L_{D_0}\cap \L=\bigcup_{\bar{\bx}\in B} \L_\bx.$ Indeed, if $L$ is a line in \L\ which meets $\bar{D}_0$, then $L\cap \bar{D}_0$ consists of a single point $\bar{\bx}\in B,$ and $L\in \L_\bx.$ Conversely, if $L\in \L_\bx$ for $\bar{\bx}\in B$, then $\bar{\bx}\in L\cap \bar{D}_0$, which implies $L\in \L_{D_0}\cap \L$. Since $\L_\bx$ is equivalent to $W_\sigma$ for each $\bar{\bx}\in B,$ the assertion follows.
\end{proof}
The subspace $D_0$ is required to meet $|\L|$ in $\deg(|\L|)$ points of depth equal to $\depth(\L).$ Since $R$ is algebraically closed, such a subspace exists.

The left-hand-side of the formula in Theorem \ref{degree} can be computed using the Pieri rule, provided $\L=\bigcap_{D\in\D} \L_D$ is a proper, generically transverse intersection. In this case each component of $\L$ has multiplicity one, so $\L=\prod_{D\in \D} W_{(c(D),0)}$ in the Chow ring of $\G(2,k)$.
\begin{corollary} Suppose $\L=\bigcap_{D\in\D} \L_D$ is a proper, generically transverse intersection, with $|\L|$ irreducible. Then the degree of $|\L|$ is determined by 
$$W_{(c(D_0),0)}\cdot \prod_{D\in \D} W_{(c(D),0)}=(\deg |\L|) \,W_\sigma,$$ where $\sigma=(k-2,k-1-\depth(\L)).$
\label{transdegree}
\end{corollary}

The intersection $\L=\bigcap_{D\in\D} \L_D$ is proper if and only if the codimension  of   \L\ in $\G(2,k)$ is equal to $\sum_{D\in \D} c(D).$ If the intersection is proper but not generically transverse,   our naive degree calculation may presumably be sharpened to take account of multiplicities, resulting in analogues of \ref{degree} and \ref{transdegree}. We leave the precise formulation and proof to the experts.

Irreducibility of $|\L(\D)|$ is a much more delicate issue. We have the following observation \cite{Pottman}.
\begin{theorem} If $\L\subseteq \G(2,k)$ is irreducible then $|\L|\subset \P^{k-1}$ is also irreducible.
\end{theorem}

\begin{proof} If $|\L|=U_1\cup U_2$, then $\L=(\L\cap {\mathcal F}_1(U_1))\cup (\L\cap {\mathcal F}_1(U_2))$, where ${\mathcal F}_1(U_i)$ is the Fano variety of lines on $U_i$.
\end{proof}

It is quite possible for $\L(\D)$ to be reducible with $|\L(\D)|$ irreducible, for instance, if $\L(\D)$ consists of the two rulings of the quadric surface in $\P^3,$ discussed below.  We have no examples of combinatorial resonance  components $\pV(\Gamma,R)$ which are not irreducible.

We close this section by briefly describing the canonical example of a projective line complex and its degenerations, to illustrate the dimension and degree formulas above. It would be interesting to see combinatorial resonance components exhiiting these phenomena. If \D\ is an arrangement of nontrivial subspaces in a vector space $K$, we denote by $\bar{\D}$ the corresponding arrangement of projective subspaces in $\bar{K}=\P(K)$. 

\subsection*{Lines in $\boldsymbol{\P^3}$} Suppose $\D=\{D_1,D_2,D_3\}$ consists of three planes in $R^4,$ so that $\bar{\D}$ consists of three lines in $\P^3.$ If $\D$ is in general position, then the three lines of $\bar{\D}$ are contained in a unique quadric surface $\Sigma$, and belong to one of the two rulings of $\Sigma.$ The line complex $\L=\L(\D)$ is the other ruling, and $|\L|=\Sigma$ is irreducible and has degree two. Every point of $|\L|$ has depth 1. The dimension of \L\ is one. This line complex is known as a {\em regulus}; its generators $D_1, D_2,$ and $D_3$ are called  {\em directrices} \cite{Veblen-Young,Pottman}. 

Suppose two lines of $\bar{\D}$ meet, say $\bar{D}_1\cap \bar{D}_2=\{\bar{\bl}\}, $ with $\bar{\bl}\not\in \bar{D}_3$. Then the plane $D_1\ast D_2$ meets $\bar{D}_3$ in a point $\bar{\be}$. In this case \L\ consists of the lines in the plane $\bl\ast D_3$ through $\bar{\bl}$ and the lines in $D_1\ast D_2$ through $\bar{\be}$. The carrier $|\L|$ is the union of the two planes $(D_1\ast D_2)$ and $(\bl\ast D_3),$ reducible, still of degree two. The two points $\bar{\bl}$ and $\bar{\be}$ have depth two; all other points of $|\L|$ have depth one. Again,  $\dim \L=1.$

If all three lines of $\bar{\D}$ meet, say at $\bar{\bl}$, but are not coplanar, then \L\ consists of all the lines in $\P^3$ through $\bar{\bl}$, $|\L|=\P^3,$  $\bar{\bl}$ has depth  three, and all other points of $|\L|$ have depth one. The dimension of \L\ is two. If the lines of $\bar{\D}$ are coplanar, i.e., $\bar{D}_3\subset D_1\ast D_2,$ then \L\ consists of all the lines in $D_1\ast D_2$, $|\L|=D_1\ast D_2,$ with every point of depth two, and $\dim \L=2$. In both of these cases $|\L|$ has degree one.

Now consider an arrangement $\D=\{D_1,D_2,D_3,D_4\}$ of four planes in general position in $R^4$. If $\bar{D}_4$ is transverse to the quadric $\Sigma=|\L(\{D_1,D_2,D_3\})|,$ then \L\ consists of the two lines of $\L(\{D_1,D_2,D_3\})$ passing through the points of $\bar{D}_4\cap \Sigma,$ and $|\L|$ is again reducible, of degree two. If $\bar{D}_4$ is tangent to $\Sigma,$ then \L\ consists of  one line, with multiplicity two. In particular, the polymatroid of \D, which tabulates the dimensions of sums of subsets of \D, is not sufficient to determine $|\L(\D)|$ up to projective equivalence.

\begin{remark} Suppose  \sG\ has rank two and $R$ is an arbitrary commutative ring. The equation for resonance derived in the last section can be interpreted in terms of line complexes. Indeed, $\l_{[n]}\be=\e_{[n]}\bl$ if and only if $$\sum_{i=1}^n \begin{vmatrix}\l_i&\e_i\\ \l_j&\e_j\end{vmatrix}=0$$ for $1\leq j\leq n.$ These are linear equations in the line coordinates of $\bl\ast\be.$ If $R$ is a field the usual methods of Schubert calculus can be used to show that these equations describe the complex of lines in the hyperplane $\x_{[n]}=0.$  As we saw in Section \ref{decomp}, this is not the case in general.
\label{ranktwocomplex}
\end{remark}

\end{section}
\begin{section}{Examples}
\label{sect-examples}

In this section we apply these ideas to several examples. The fine structure revealed in Section \ref{sect-lines} is not apparent in resonance varieties over \C. In this case $|\L(\D_\Gamma)|$ is known to be linear, and  $\L(\D_\Gamma)$ is the complex of all lines in $|\L(\D_\Gamma)|$. Resonance components supported on complexified real arrangements are lines. Indeed there is only one known example, supported on an arrangement of rank greater than two, for which $\L(\D_\Gamma)$ consists of more than a single line. Nontrivial line structure emerges over fields of positive characteristic. In particular, we will see that the Hessian arrangement supports a resonance component over $\bar{\Z}_3$ that is a cubic threefold, with interesting line structure. 

Throughout this section, we assume $R$ is a field, unless otherwise specified. To begin, we recall the various special features of resonance varieties over fields of characteristic zero. Let \sG\ be a simple matroid on ground set $[n]$. Let $\Gamma$ be a graph with vertex set $[n]$. Recall from Section \ref{decomp},  $\X_\Gamma(\sG)$ is the set of nontrivial lines of $\sG$ which are not cliques of $\Gamma$, and $K(\Gamma,R)=\{\bx\in R^n \ | \ \x_X=0 \ \text{for every} \ X\in \X_\Gamma(\sG)\}.$ Let $\Delta=\{\bx\in R^n \ | \ \x_{[n]}=0\}$ and $K_0(\Gamma,R)=K(\Gamma,R) \ \cap \ \Delta$.

\begin{theorem} Suppose $R$ is a field of characteristic zero, and $K_0(\Gamma,R)$ has dimension at least two. Then
\begin{enumerate}
\item $V^1(\Gamma,R)=K_0(\Gamma,R),$ a linear subspace of $R^n.$
\item If $\bl\in V^1(\Gamma,R)$, then $Z(\bl)=V^1(\Gamma,R)=V^1_d(\Gamma,R),$ for $d=\dim Z(\bl)-1.$
\item The dimension of $V^1(\Gamma,R)$ is one less than the number of blocks of $\Gamma$. In particular $\Gamma$ has at least three blocks.
\item  If $V^1(\Gamma,R)\not = V^1(\Gamma',R)$, then $V^1(\Gamma,R)\cap V^1(\Gamma',R)=0.$
\item If $X\in \X_\Gamma(\sG)$, then every block of $\Gamma$ meets $X$.
\item If $\bl\in V^1(\Gamma,R)$ with  $\supp(\Gamma)=[n]$, then $\bl$ is constant on blocks of $\Gamma.$
\item If $\sG$ is realizable over $\Re$, then $\Gamma$ has at most three blocks.
\end{enumerate}
\label{yuz}
\end{theorem}
The first six results are proved in \cite{LY}, see also \cite{Y9}. The last is  a consequence of (v), as shown in \cite{CorFor}. Note that the hypothesis on $\Gamma$ really depends only on $\X_\Gamma(\sG).$ In particular $\Gamma$ is not assumed to be neighborly - this follows from (v). In \cite {LY} these graphs arise from a block  (direct sum) decomposition of the matrix $Q=I^TI-J,$ where $I$ is the incidence matrix defining $K(\Gamma,R)$ and $J$ is the matrix of all one's. $Q|_\Delta$ is positive-definite, and $\ker Q|_\Delta=\ker I|_\Delta,$ when $R$ has characteristic zero. 

\medskip
We proceed with a few elementary observations, with proofs left to the reader. These results treat the trivial cases, which encompass almost all known examples. Suppose \D\ is a subspace arrangement in a  vector space $K$ of dimension $k>0.$  We define the {\em proper part} $\D_0$ of \D\ by $\D_0=\{D\in \D \ | \ \codim_K(D)>1\}.$ 

\begin{theorem} \begin{enumerate}
\item $\L(\D)=\L(\D_0)$
\item If $|\D_0|=\emptyset$, then  $|\L(\D)|=\bar{K}$, with every point of depth $(k-1)$.
\item If $\D_0=\{D\}$, then $|\L(\D)_\bl|=\bl\ast D$ for every $\bl\in K$, and  $|\L(\D)|=\bar{K}.$ Every point of $\bar{K}$ has depth $\dim(D)$.
\item If $\D_0=\{D_1,D_2\}$ with $D_1\not = D_2$, then $|\L(\D)|=D_1\ast D_2.$ Points of $\bar{D}_1\cap\bar{D}_2$ have depth $\dim (D_1 + D_2)$; all other points have depth $\dim (D_1\cap D_2) +1$.
\end{enumerate}
\end{theorem}

In general $\D$ may contain one or more 1-dimensional subspaces. We will call such elements the {\em poles} of \D. 

\begin{theorem} Let $\D=\D_\Gamma$ for $\Gamma\in \NP(\sG,R)$. Suppose $\D$ contains a pole $R\bl$. Let $\L=\L(\D)$. Then 
\begin{enumerate} 
\item $\L=\L_\bl$.
\item $\pV(\Gamma)=|\L_\bl|=Z_\Gamma(\bl)=\bigcap_{D\in\D_0} \bl\ast D$ 
\item $\pV_d(\Gamma,R)=0$ for $d>\depth(\bar{\bl})$.
\item If $\depth(\bar{\bl})\geq 2$ then $\pV_d(\Gamma,R)=\{\bar{\bl}\}$ for $2\leq d\leq \depth(\bar{\bl})$.
\item If $\D_0=\{R\bl\}$ then $\pV(\Gamma)=\bar{K}$.
\end{enumerate}
\label{single}
\end{theorem}

\begin{corollary} Let $\D=\D_\Gamma$ for $\Gamma\in \NP(\sG,R)$. Suppose $\D$ has poles $R\bl_1, \ldots, R\bl_n,$ with $n\geq 2.$ Then $|\L(\D)|=\emptyset$ unless $\bar{\bl}_1, \ldots, \bar{\bl}_n$ are collinear, in which case 
\begin{enumerate}
\item $\pV(\Gamma)=\bl_1\ast \bl_2$ and
\item $V^1_d(\Gamma,R)=0$ for $d>1$.
\end{enumerate}
\end{corollary}

By these results, if the arrangement of directrices $\D_\Gamma$ has fewer than three subspaces of codimension greater than one, or contains an element of dimension one, then the combinatorial component $V(\Gamma,R)$ is linear. 

\subsection*{Local components in $\mathbf \R^1(\sG,R)$} Suppose \A\ is a pencil of $n\geq 3$ lines in the plane. Then $\sG(\A)$ is an $n$-point line.
Let $\Gamma=\emptyset$. Then $\X_\Gamma(\sG)$ consists of the single (nontrivial) line in $\sG$, and the point-line incidence matrix has rank 1. Then $K=K(\Gamma,R)$ has dimension $n-1$. The directrices $D_{\{i\}}$ are all hyperplanes, so $\pV(\Gamma,R)=|\L_{\D_\Gamma}|=\bar{K}$. 

For arbitrary $\sG$, a resonant weight supported on a flat of rank two is called ``local." For $X\in \X_0(\sG),$ define the graph $\Gamma_X$ by $$\Gamma_X=\{\{i,j\} \ | \ |\{i,j\}\cap X|\leq 1\}.$$ Then every point of $[n]-X$ is a cone vertex, $\X_\Gamma(\sG)=\{X\}$, and $V^1(\Gamma,R)=K(\Gamma,R)$ is a linear subspace of dimension $|X|-1$ as above. Suppose $R$ is algebraically closed. Then $V^1(\Gamma,R)$ is  irreducible.  In fact, the proof given in \cite{F7} can be adapted to show that $V^1(\Gamma,R)$ is an irreducible component of $\R^1(\sG,R)$ in this case. These are the {\em local components} of $\R^1(\sG,R)$.

\medskip
In the examples below, we give arrangements in terms of their defining polynomials. We order the hyperplanes according to the order of factors in the defining polynomial. We will use $\alpha, \beta, \ldots$ in place of 2-digit labels. We illustrate some of the examples using affine matroid diagrams; the interpretation should be clear, but the reader may consult \cite{Ox} for a detailed explanation. We specify graphs by listing their maximal cliques in block notation. We will write vectors over $\Z_2$ as bit strings.

In almost all known examples, $K=K(\Gamma,R)$ has dimension two. Then $\D_0(\Gamma)$ is empty, and $\pV(\Gamma,R)=\bar{K}$ is a line. In particular this is always the case for non-local resonance in real arrangements (or real-realizable matroids) over fields of characteristic zero: by Theorem \ref{yuz} (vii), $\Gamma$ has three blocks, and then by (iii), $K(\Gamma,R)$ has dimension two.

We exhibit one such example for future reference, the canonical example of non-local resonance \cite{F7}.

\begin{example} The braid arrangement in $\Re^3$ has defining polynomial $$Q(x,y,z)=(x+y)(x-y)(x+z)(x-z)(y+z)(y-z),$$ with underlying matroid $\sG$ isomorphic to the cycle matroid $K_4$ of the complete graph on four vertices. The nontrivial lines in $\sG$ are $136, \ 145, \ 235, \ \text{and} \ 246.$ A non-local resonance component arises from the neighborly partition $\Gamma=12|34|56.$ 
\begin{figure}
\end{figure}
We have $\X_\Gamma(\sG)=\X_0(\sG)$, and the $4\times 6$ point-line incidence matrix has rank 4, over any field $R$. Thus $K=K(\Gamma,R)$ has dimension two, so $\pV(\Gamma,R)=\bar{K}$.  The arrangement of directrices $\bar{\D}_\Gamma$ consists of the three collinear points $\bar{\bl}_1, \bar{\bl}_2, \bar{\bl}_3$, where $\bar{\bl}_1=(1,1,0,0,-1,-1), \bar{\bl}_2=(0,0,1,1,-1,-1),$ and $\bar{\bl}_3=(1,1,-1,-1,0,0)$.  In particular $(\bl_i,\bl_j)$ is a resonant pair for $i\not =j$. It is useful to point out here  that, in case $\charac(R)=2,$ each $\bl_i$ is a sum of characteristic functions of blocks of $\Gamma$: $\bl_1=110011, \bl_2=001111,$ and $\bl_3=111100.$ 
\label{braid}
\end{example}

The following example, discovered by C. Olive and E. Samansky, illustrates that the arrangement of directrices can be in special position in $K$, so that the corresponding line complexes do not intersect properly. (If $K$ is replaced by the span of the directrices, the intersection becomes proper.) This example shows that the vanishing of $(a_\bl\we a_\be)_X$ for all nontrivial lines $X$ does not imply $a_\bl\we a_\be=0,$ as is the case in characteristic zero.

\begin{example} Consider the arrangement of 10 planes in $\Re^3$ with defining polynomial 
$$
(x+z)(x-z)(x+y)(x-y)(2x+z)(2x-z)(x+2y)(x-2y)(2y+z)(2y-z).$$
Then $\X_0(\sG)=\{17\alpha,189,279,28\alpha,35\alpha,369,459, 46\alpha,1256,3478\}.$ Let  $$\Gamma=1234|5678|9\alpha.$$ Then $\Gamma$ is neighborly. The hypothesis  of Theorem \ref{yuz} is not satisfied; indeed, $\dim K(\Gamma,R)=1$ if $R$ is a field of characteristic zero. Suppose $R$ is a field of characteristic two. Then $K=K(\Gamma,R)\cong R^4$. In $\bar{K}\cong \P^3$ the arrangement of directrices $\bar{\D}=\bar{\D}_0$ consists of two poles  $\bl_1=\bar{1111000011}$ and $\bl_2=\bar{0000111111}$, and the line $L=\bx\ast\bn$, where $\bx=1100001100$ and $\bn=0011110000.$ Observe that the line $\bl_1\ast\bl_2$ meets $L,$ at the point $R\, 1111111100. $ Thus $\pV(\Gamma,R)$ is the line $\bl_1\ast\bl_2.$ Note that $V^1(\Gamma,R)\not = K_0(\Gamma,R)$, contrary to Theorem \ref{yuz}(i).
\end{example}

The next example, which  exhibits higher order, non-local resonance that only appears in characteristic two, was found by D.~Matei and A.~Suciu \cite{MatSuc1}. This example provided the original motivation for the present study. It illustrates that resonance in characteristic two is governed by incidences among submatroids of $\sG$, a phenomenon that has farther-reaching consequences in Example \ref{deleted1} below.

\begin{example} Consider the real arrangement defined by $$Q(x,y,z)=(x+y)(x-y)(x+z)(x-z)(y+z)(y-z)z.$$ Its underlying matroid $\sG$ is the non-Fano plane, with nontrivial lines $$136, \ 145, \ 235, \ 246, \ 347, \ 567.$$ If $R$ is a field of characteristic zero, no transitive neighborly graph $\Gamma$ with $\supp(\Gamma)=[n]$  satisfies $\dim K(\Gamma,R)>1$. (In fact, \sG\ has no neighborly partitions with only three blocks, as would be required by Theorem~\ref{yuz}(vii).) Suppose $R$ is a field of characteristic two and $\Gamma=127|3|4|5|6$. Then $\Gamma$ is neighborly, and $\X_\Gamma(\sG)=\X_0(\sG)$. The $6\times 7$  incidence matrix for $\X_\Gamma(\sG)$ has rank 4 over $R$. Hence $K= K(\Gamma,R)$ has dimension 3, and $\bar{K}$ is a plane. 
The directrices corresponding to the singleton blocks of $\Gamma$ are lines in this plane, while $D_{127}$ is a pole $\bar{\bl},$ with $\bl=0011110.$ Thus $\D_0=\{D_{127}\}$, and $\pV(\Gamma,R)=\bar{K}$ by Theorem \ref{single}. The pole $\bar{\bl}$ has depth two, and $\pV(\Gamma,R)=|\L_\bl|=\P(Z_\Gamma(\bl)).$

\begin{figure}[h]
\begin{center}
\includegraphics[height=4cm]{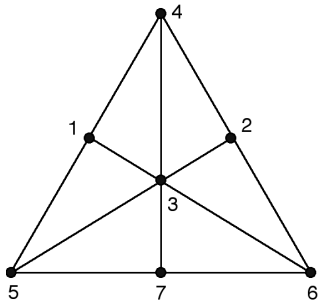}
\end{center}
\caption{The non-Fano plane}
\label{figure-nonfano}
\end{figure}

The diagram of $\sG$  appears in Figure~\ref{figure-nonfano}.  Note that deleting, respectively, points 1, 2, and 7 give submatroids of $\sG$ isomorphic to $K_4$, the matroid of Example~\ref{braid}. Comparing with that example, observe that the special weight $\bl$ is the characteristic function of the intersection of these three submatroids, and in fact is the sum of characteristic functions of two of the three blocks in each of them. For each $i=1,2,7$, there is a resonant pair $(\bl,\be_i)$ supported on $\sG-\{i\}$. Since $\be_1+\be_2+\be_7=0,$ one gets $\dim Z(\bl)=3$. Note that there are no resonant pairs $(\bl,\be)$ with $\be\in (\Z_2)^7$ which are supported on \sG. So $\Gamma=\Gamma_\bl$ does not coincide with any $\Gamma_{(\bl,\be)}$ for $\be\in (\Z_2)^7$ - it is necessary to pass to a field extension (e.g., $\bar{\Z}_2$) to find a generic partner for $\bl.$

The incidence geometry that yields resonance over $\Z_2$ is also reflected in the characteristic varieties $\V_k^1(\A)$: the three components of $\R^1(\A,\C)$ corresponding to the $K_4$ submatroids exponentiate to three 2-tori in $\V_1^1(\A)\subseteq (\C^*)^7,$ which intersect at the point $(1,1,-1,-1,-1,-1,1)=\exp(2\pi i\bl/2).$ This point is precisely $\V^1_2(\A)$.  This was the first known example of a component of a characteristic variety which does not pass through the identity \cite{CS3}. 
\label{nonfano}
\end{example}

\subsection*{The deleted $B_3$ arrangement} In \cite{Suc} A. Suciu introduced the ``deleted $B_3$ arrangement," obtained by deleting one plane from the reflection arrangement of type $B_3$. The defining polynomial is given by 
$$Q(x,y,z)=(x+y+z)(x+y-z)(x-y-z)(x-y+z)(x-z)x(x+z)z.$$

Suciu showed that the characteristic variety $\V_1(\A)\subseteq (\C^*)^n$ has a one-dimensional component which does not pass through $(1,\ldots, 1)$. Since the components of $\V^1(\A)$ containing  $(1,\ldots, 1)$ are tangent to the components of $\R^1(\A,\C)$, they all have dimension at least two, by Theorem \ref{yuz}. To that point, no arrangements had been found with other than $0$-dimensional components in $\V_1(\A)$ away from $(1,\ldots,1)$. In the next example we see that the same incidence structure that gives rise to this translated component in $\V^1(\A)$ yields components of $\R^1(\A,R)$ with nontrivial intersection, for $\charac(R)=2$, in contrast to Theorem \ref{yuz} (iv).
\begin{example} Let \sG\ be the matroid of the deleted $B_3$ arrangement, illustrated in Figure \ref{db3diagram}.

\begin{figure}[h]
\begin{center}
\includegraphics[height=4cm]{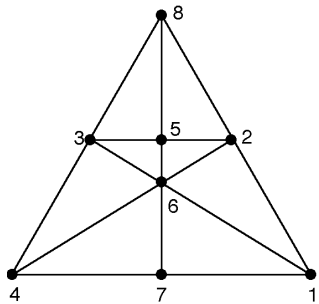}
\end{center}
\caption{The deleted $B_3$ matroid}
\label{db3diagram}
\end{figure}

The deletions $\sG_7=\sG-\{7\}$ and $\sG_5=\sG-\{5\}$ of $\sG$ are copies of the non-Fano plane, and their intersection $\sG_{57}$ is a copy of $K_4$. Let $R$ be an algebraically closed field of characteristic two. With Example \ref{nonfano} in mind, we set $\bl_1=01100101$ and $\bl_2=10010101.$ Then $Z(\bl_1)=V^1(\Gamma_1,R),$ where $\Gamma_1=1457|27|37|67|78$ corresponds to the neighborly partition $145|2|3|6|8$ of $\sG_7.$  Similarly, $Z(\bl_2)=V_1(\Gamma_2,R)$ where $\Gamma_2=2357|15|45|56|58,$ with $\supp(\Gamma_2)=\sG_5.$

Now let $\be=11110000.$ Observe that $\be$ is supported on $\sG_{57}$, and in fact is a sum of (characteristic functions of) blocks of the neighborly partition $14|23|68$ of $\sG_{57}$. Each of $\bl_1=01100101$ and $\bl_2=10010101$ is a sum of blocks of the same partition. Thus $(\bl_1,\be)$ and $(\bl_2,\be)$ are resonant pairs, as in Example \ref{braid}, and so $\be\in Z(\bl_1)\cap Z(\bl_2).$ Indeed,  $Z(\bl_1)\cap Z(\bl_2)=Z(\be).$ In particular $V_1(\Gamma_1,R)\cap V_1(\Gamma_2,R)$ is nontrivial. Since $\sG$ itself does not support any neighborly partitions, $V_1(\Gamma_i,R)$ is indeed an irreducible component of $\R^1(\sG,R)$ for $i=1,2.$

In Figure \ref{db3picture} is a picture of $\pV(\Gamma_1)\cup \pV(\Gamma_2)$, with the line structure indicated in bold.

\begin{figure}[h]
\begin{center}
\includegraphics[height=6cm]{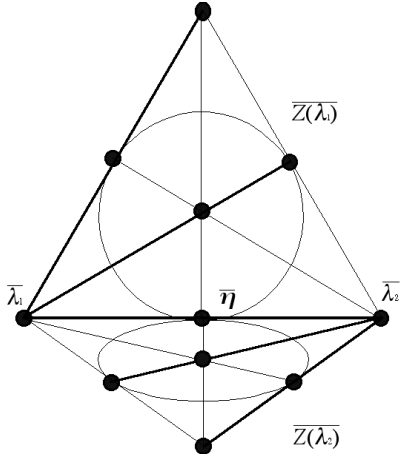}
\caption{A nontrivial intersection}
\label{db3picture}
\end{center}
\end{figure}
Note that $\bl_1\in Z(\bl_2)$ but $Z(\bl_1)\not\subseteq Z(\bl_2).$ That is, there exists $\bx\in R^7$ such that $a_{\bx}\we a_{\bl_1}=0=a_{\bl_1}\we a_{\bl_2},$ but $a_{\bx}\we a_{\bl_2}\not =0,$ suggesting the possible existence of nontrivial triple Massey products over $\Z_2.$  This phenomenon is expected when $\R_2^1(\sG,R)$ is nontrivial and strictly contained in $\R_1^1(\sG,R),$ as in Example~\ref{nonfano}. For that example, D.~Matei informs us that there are no nontrivial triple Massey products. In the present example, we have $\R_2^1(\sG,R)=0.$ We do not know whether $\langle \bx,\bl_1,\bl_2\rangle$ is nontrivial for some $\bx.$

Referring to \cite{Suc} we find that $\bl_1$ and $\bl_2$ correspond via exponentiation $\bl\mapsto \exp(2\pi i \bl/2)$ to the two points $\rho_1$ and $\rho_2$ that comprise the second characteristic variety $\V^1_2(\A)$, points which lie on the translated component $$C=\{(t,-t^{-1},-t^{-1},t,t^2,-1,t^{-2},-1) \ | \ t\in \C^*\}$$ of $\V_1(\A).$ 
The difference $\bl_1-\bl_2$ is $\be,$ which exponentiates to $\rho_1\rho_2^{-1}.$ Thus $Z(\be)$ exponentiates to the one-dimensional subgroup whose coset by $\rho_i$ is $C$.
The diagram \cite[Figure 6]{Suc} indicates that the same overlapping of non-Fano's and $K_4$'s in $\sG$ gives rise to the existence of $\be$ and to that of $C$. 
\label{deleted1}
\end{example}

None of the weights in the preceding example have full support. Indeed, the underlying matroid \sG\ has no neighborly partitions, so there can be no resonant pairs with full support over any integral domain. Yet the translated torus $C$ is not supported on any proper subarrangement. In conversations with D.~Cohen attempting to better understand this situation, we discovered a resonant pair over $\Z_4$, with full support. It was to understand this phenomenon that the theory of Section \ref{decomp} was developed  for arbitrary commutative rings. A variation on that example is presented next; many of the anticipated pathologies exhibit themselves.

\begin{example} Let $R$ be a ring of characteristic 4, say $R=\Z_4.$ Let $$\bl=(1,1,1,1,2,2,2,2)\ \text{and} \ \be=(2,3,1,0,0,1,2,3).$$ Using the Smith Normal Form of $d_\bl: A^1 \to A^2$ one finds that $Z(\bl)=\ker(d_\bl)$ is a free $R$-module of rank two, with basis $\{\bl,\be\}$.  Thus $(\bl,\be)$ is a resonant pair, with $\supp(\bl,\be)=[8].$
The graph $\Gamma=\Gamma_{(\bl,\be)}$ has edges $15,27,37,45,57,68.$ So $\Gamma$ is not transitive, i.e., does not yield a partition of $[8].$ The blocks of $\Gamma$ all have size two, with $5$ and $7$ each belonging to three blocks. Referring to Figure \ref{db3diagram}, we see that $\Gamma$ is indeed a neighborly graph. In particular, $\bl_X$ is parallel to $\be_X$ for every trivial line $X$. We have $\X_\Gamma(\sG)=\X_0(\sG)=\{128, 136, 147, 235, 246, 348,5678\}$. Then one checks that $\l_X=0=\e_X$ (modulo 4) for every $X\in\X_0(\sG)$, except $X=5678.$ In fact $\e_{5678}=2.$ But  it is true that $\l_{5678}\be_{5678}=\e_{5678}\bl_{5678}$ as required by Theorem \ref{res-pair}.

Let $\be'=(1,3,3,1,2,0,2,0)=\bl+2\be\in Z(\bl).$ Then $\be'$ is not parallel to \bl, so $(\bl,\be')$ is also a resonant pair. This in spite of the fact that \bl\ and $\be'$ are linearly dependent: $2\bl+2\be'=0.$  The graph $\Gamma_{(\bl,\be')}$ is the union of three complete graphs, on vertices 1457, 2357, and 5678, and $\bl\in V^1(\Gamma',R)$.   Again $\bl_X$ is parallel to $\be'_X$ for every trivial line $X$, and now $\l_X=0=\e'_X$ for every $X\in\X_0(\sG).$ This is the pair we found with Cohen. Vertices $5$ and $7$ are cone vertices which lie in $\supp(\bl,\be'),$ so again $\Gamma'$ is not transitive. In fact $\Gamma'$ is the graph of the combinatorial component associated with the copy $\sG_{57}$ of $K_4$, and induces a neighborly partition on $[8]-\{5,7\}.$. All the resonant weights supported on $\sG_{57}$ lie in $V^1(\Gamma),$ though we can not tell whether $V^1(\Gamma')\subset V^1(\Gamma).$

Referring again to the translated component $C$ in $\V_1(\A)$, observe that \bl\ exponentiates to a point $\exp(2\pi i \bl/4)=i^\bl$ that lies on $C,$ while $\exp(2\pi i\be'/4)=i^{\be'}$ generates the subgroup of $(\C^*)^n$ corresponding to $C$. On the other hand $\exp(2\pi i\be/4)=i^{\be}$ itself seems to have no relation to $C$.
\label{db3-z4}
\end{example}

\subsection*{The Hessian arrangement} Finally we present an example having nonlinear components in $\R^1(\sG,R).$ Let \A\ be the Hessian arrangement in $\C^3,$ corresponding to the set of twelve lines passing through the nine inflection points of a nonsingular cubic in $\P^2(\C)$ \cite[Example 6.30]{OT}.  The underlying matroid is the deletion of one point from $PG(2,3),$ the projective plane over $\Z_3.$
We choose a labelling so that 

$$\X_0(\sG)=\{148\gamma, 159\alpha, 167\beta, 247\alpha,258\beta, 269\gamma, 349\beta, 357\gamma, 368\alpha\}.$$

\begin{example} Let $R=\bar{\Z}_3, $  and let $\Gamma=123| 456 | 789|\alpha\beta\gamma$. Then we have $\X_\Gamma(\sG)=\X_0(\sG).$ The $9 \times 12$ point-line incidence matrix has rank six over $R$, so $\dim K(\Gamma,R)=6.$ For each block $S$ of $\Gamma,$ the corresponding directrix has dimension three. Then the projective arrangement of directrices $\bar{\D}_\Gamma$ consists of four planes in $\P^5.$ (For $R=\C$ the  directrices are four lines in $\P^2.$)

The placement of these four planes is special: each meets the other three in three collinear points. It follows that the six points of intersection are coplanar, and are the six points of intersection of four lines in general position in that plane. The join of the four planes in $\bar{\D}_\Gamma$ is a $\P^4$ in $\P^5.$ 

A {\em Macaulay2} computation tells us that the associated ruled variety is an irreducible cubic hypersurface in $\P^4.$ This can be confirmed by analyzing the line complex $\L=\L(\D_\Gamma)=\bigcap_{P\in \D_\Gamma} \L_P$ using the Schubert calculus in $\P^4.$ 

By Theorem \ref{cell}, a plane $\bar{P}$ in $\P^4$ determines a line complex $\L_P$ equivalent to the Schubert variety $W_{(1,0)}$. If the intersection $\L=\L(\D_\Gamma)=\bigcap_{P\in \D_\Gamma} \L_P$ is proper, then $\L$ is rationally equivalent to 
$W_{(1,0)}^4$, which equals $3W_{(3,1)}+2W_{(2,2)}$ by the Pieri rule.

We can show the intersection is proper  by the following {\em ad hoc} argument. The codimension of $\L$ in $G(2,5)$ is at most 4, since that is the codimension of $W_{(1,0)}^4,$  and codimension does not increase under degeneration. One sees without much difficulty that $\L$ has depth one. Then, $\codim \L=\codim |\L|+3$ by Theorem \ref{dimensions}. Since $|\L|$ is easily seen to be a proper subvariety of $\P^4$ we conclude $\codim \L=4,$ as desired. 

We have no method to show the intersection is generically transverse.
Assuming it is, the degree of $|\L|$ is determined as follows. We have shown  $\dim |\L|=3.$ So we take a subspace $D_0$ of dimension two, corresponding to a line in $\P^4$, and calculate the intersection of $\L_{D_0}\cong W_{(2,0)}$ with $\L$:
$$\L\cap \L_K\sim (3W_{(3,1)}+2W_{(2,2)})\cdot W_{(2,0)}=3W_{(3,3)}.$$
Then $\deg |\L|=3$ by Theorem \ref{degree}. It would be nice to complete the argument by establishing the transversality, or by generalizing Corollary \ref{transdegree}, but this is beyond our ken at the present.

We have not found a direct argument, {\em ad hoc} or otherwise, to show $|\L|$ is irreducible.

H.~Schenck analyzed the cubic threefold $|\L|$ using {\em Macaulay2}. The plane containing the six intersection points of the directrices is singular in $|\L|$. These are the points of depth two. The quadric in that plane consisting of the four lines containing the intersection points is an embedded component. This variety is apparently  related to  more familiar threefolds \cite{Harris,Todd}. As it seems to be a primal geometric object, at least in the characteristic-three universe, it deserves a more precise description.  Again this undertaking is left for the experts.
\label{hessian}
\end{example}

We close by listing a few problems.

\begin{problem} Find a matroid \sG\ with neighborly partition $\Gamma$ such that the corresponding line complex $\L(\D_\Gamma)$ is a regulus, i.e., so that the directrices are three skew lines in $\P^3.$
\end{problem}

\begin{problem} Find a matroid \sG\  with neighborly partition $\Gamma$ such that the corresponding combinatorial component $\pV(\Gamma,R)$ is reducible.
\end{problem}

\begin{problem} Determine the elementary divisors of the point-line incidence matrix of an arbitrary rank-three matroid in terms of more familiar invariants.
\end{problem}

\begin{problem} Explain how resonance in characteristic $N$ gives rise to $N$-torsion points and/or translated components in the characteristic variety.
\end{problem}

\begin{problem} Show how the linearity and trivial intersection properties of resonance varieties over fields of characteristic zero can be deduced from the neighborly graph/line complex description of resonance components.
\end{problem}

\begin{problem} Determine the codimension and degree of $\L(\D)$ for an arbitrary arrangement \D\ of projective subspaces.
\end{problem}

\end{section}

\begin{ack} This research was helped along by discussions with many people. In the beginning, work by my REU students Cahmlo Olive and Eric Samansky, in the summer of 2000, provided interesting examples to scrutinize. I am grateful to Sergey Yuzvinsky, Hiroaki Terao, Dan Cohen, Alex Suciu, Hal Schenck, and the other participants in two mini-workshops at Oberwolfach (March, 2002 and November, 2003) for their help. Frank Sottile helped me to gain a rudimentary understanding of intersection theory.
\end{ack}

\providecommand{\bysame}{\leavevmode\hbox to3em{\hrulefill}\thinspace}
\providecommand{\MR}{\relax\ifhmode\unskip\space\fi MR }
\providecommand{\MRhref}[2]{%
  \href{http://www.ams.org/mathscinet-getitem?mr=#1}{#2}
}
\providecommand{\href}[2]{#2}

\end{document}